\magnification 1200
\documentstyle{amsppt}
\NoBlackBoxes

\topmatter

\loadeufm


\define \im{\text{\rm Im }}
\define \re{\text{\rm Re }}

\define \bR{\Bbb R}
\define \bZ{\Bbb Z}
\define \bC{\Bbb C}
\define \scrM{\Cal M}
\define \scrA{\Cal A}
\define \scrO{\Cal O}
\define \scrV{\Cal V}
\define \scrL{\Cal L}
\define \ssneq{\subsetneqq}
\define \hol{\text{\rm Hol}}
\define \alg{\text{\rm Alg}}
\def\cnn {\bC^N}
\def\ctnn {\bC^{2N}}
\def\cn {\bC^n}
\def \po {p_0}
\def\r {\rho}
\def\pa {\partial}
\def\a {\alpha}

\def\z {\zeta}
\def\d {\delta}
\def\t {\tau}
\def\ch {\chi}
\def\cd {\bC^d}

\def\J {\text {\rm Jac\ }}

\def\areg {A_{\text {\rm reg}}}
\def\acr {A_{\text {\rm CR}}}
\def\mcr {M_{\text {\rm CR}}}
\def\dim {\text {\rm dim}}

\title Algebraicity of holomorphic mappings
between real algebraic sets in
${\bold C}^n$ \endtitle
\rightheadtext{Algebraicity of holomorphic
mappings}
\leftheadtext{M.~S.~Baouendi, P.~Ebenfelt, and
L.~P.~Rothschild}
\author M. S.
Baouendi\footnote {Partially supported by National
Science Foundation Grant DMS
95-01516\hfill\break},  P.
Ebenfelt\footnote{Supported by a grant from the
Swedish Natural Science Research
Council.\newline}, and Linda Preiss Rothschild$^1$
\endauthor
\address Department of Mathematics-0112,
University of California, San Diego,\hfill\break La
Jolla, CA 92093-0112\endaddress
\email sbaouendi\@ucsd.edu,
pebenfel\@math.ucsd.edu,
 lrothschild\@ucsd.edu\endemail

\toc
\widestnumber\head{\S 3.3.}
\head \S 0. Introduction\endhead
\head \S 1. Holomorphic nondegeneracy of
real analytic manifolds\endhead
\subhead \S1.1 Preliminaries on real
submanifolds of $\cnn$
\endsubhead
\subhead \S1.2 Holomorphic nondegeneracy and
its propagation 
\endsubhead
\subhead \S 1.3 The Levi number and essential
finiteness 
\endsubhead
\subhead \S 1.4 Holomorphic
nondegeneracy of real algebraic sets
\endsubhead
\head \S 2. The Segre sets of a real analytic CR
submanifold \endhead
\subhead \S2.1. Complexification of $M$,
involution, and projections\endsubhead
\subhead \S2.2. Definition of the Segre sets of $M$
at
$p_0$\endsubhead
\subhead \S2.3. Homogeneous submanifolds of
CR dimension 1\endsubhead
\subhead \S2.4. Homogeneous submanifolds of
arbitrary CR dimension\endsubhead
\subhead \S2.5. Proof of Theorem
2.2.1\endsubhead

\head \S 3. Algebraic properties of
holomorphic mappings between real algebraic
sets\endhead
\subhead \S3.1. A generalization of
Theorems 1 and 4\endsubhead
\subhead \S3.2. Propagation of
algebraicity\endsubhead
\subhead \S3.3. Proof of Theorem
3.1.2\endsubhead
\subhead \S3.4. Proof of Theorem
3.1.8\endsubhead
\subhead \S3.5. An example\endsubhead
\subhead \S3.6.  Proofs of
Theorems 1 through 4\endsubhead
\endtoc

\abstract Let $A\subset \bC^N$ be an irreducible
real algebraic set.  Assume that there exists
$\po \in A$ such that $A$ is a minimal,
generic, holomorphically nondegenerate submanifold
at
$\po$.  We show here that if
$H$ is a germ at
$p_1 \in A$ of a holomorphic mapping from $\bC^N$
into itself, with Jacobian $H$ not identically $0$,
and
$H(A)$ contained in a real algebraic set of the same
dimension as $A$, then $H$ must extend to all of
$\bC^N$ (minus a complex algebraic set) as an
algebraic mapping. Conversely, we show that
for any ``model case" (i.e., $A$
given by quasi-homogeneous real polynomials), the
conditions on $A$ are actually necessary for the
conclusion to hold.
\endabstract

\endtopmatter

\document 
\heading 0. Introduction\endheading

A subset $A\subset \cnn$ is a {\it real algebraic
set} if it is  defined by the vanishing of real
valued polynomials in
$2N$ real variables. By $\areg$ we mean the
regular points of $A$, i.e. the points at
which $A$ is a real submanifold of $\cnn$. If
$A$ is irreducible, we write $\dim\ A =
\dim_\bR\ A$ for the dimension of the real
submanifold $\areg$. A germ of a holomorphic
function $f$ at a point
$\po \in \cnn$ is called {\it algebraic} if
it satisfies a polynomial equation of the
form
$$a_K(Z)f^K(Z)+\ldots+a_1(Z)f(Z)+ a_0(Z)
\equiv 0, 
$$ where the $a_j(Z)$ are holomorphic
polynomials in $N$ complex variables with
$a_K(Z) \not\equiv 0$.  A real analytic
submanifold in $\cnn$ is called {\it
holomorphically degenerate} at $p_0 \in M$
if there
  exists  a  germ at
$\po$ of a holomorphic vector field, with
holomorphic coefficients,  tangent
 to $M$ near $\po$, but not vanishing
identically on $M$ ; otherwise, we say that
$M$ is {\it holomorphically nondegenerate} at
$\po$ (see
\S 1). In this paper, we shall give
conditions under which a germ of a
holomorphic map in
$\cnn$, mapping an irreducible real
algebraic set
$A$ into another of the same dimension, is
actually algebraic. We shall now describe
our main results.

\proclaim {Theorem 1} Let $A \subset
\cnn$ be an irreducible real algebraic set. 
Suppose the following two conditions hold.
\roster
\item  $A$ is holomorphically nondegenerate
at every point of some nonempty relatively
open subset of
$A_{\text {\rm reg}}$.
\item  If $f$ is a germ, at a point in
$A$, of a  holomorphic algebraic function in
$\cnn$  such that the restriction of
$f$ to $A$
  is real valued, then $f$ is constant.  
\endroster Then if 
$H$ is a holomorphic map from an open
neighborhood in
$\cnn$ of a point
$\po \in A$ into $\cnn$, with
$\J H|_A
\not\equiv 0$, and mapping $A$ into another real
algebraic set $A'$ with {\rm dim}
$A'$ = {\rm dim} $A$, necessarily the map
$H$ is algebraic.
\endproclaim 

We shall show that the conditions (1) and
(2) of Theorem  1  are essentially necessary
by giving a converse to Theorem 1.  For
this, we need the  following definitions. If
$M$ is a real submanifold of $\cnn $ and
$p
\in M$, let
$T_pM$ be its real tangent space at
$p$, and let
$J$ denote the anti-involution of the
standard complex structure of
$\cnn$.  We say that $M$ is {\it CR} (for
Cauchy-Riemann) at $p$ if
$\dim_\bR (T_qM +JT_qM)$ is constant for $q$
in a neighborhood of $p$ in
$M$.  If $M$ is CR at $p$, then
$\dim_\bR T_pM\cap JT_pM = 2n$  is even and
$n$ is called the {\it CR dimension} of $M$
at
$p$.  We shall say that an  algebraic
manifold
$M\subset \cnn$  is  {\it homogeneous } if
it is  given by the vanishing of
$N-\dim\ M$ real valued polynomials, whose
differentials are linearly independent at
$0$, and  which are homogeneous with respect
to some set of weights (see \S 3.6).  

\proclaim {Theorem 2}  Let $A
\subset
\cnn$ be an irreducible real algebraic set,
and let {\rm (1)} and {\rm (2)} be the
conditions of Theorem \text {\rm 1}.
Consider the following property.
\roster
\item "(3)" For every $\po \in
\areg$ at which
$A$ is CR there exists a germ of a
nonalgebraic biholomorphism $H$ of
$\cnn$ at
$\po$ mapping
$A$ into itself  with
$H(\po) = \po$.
\endroster  If {\rm (1)} does not hold then
{\rm (3)} holds. If \text {\rm (1)}  holds,
but \text {\rm (2)} does not hold, let
$f$ be a nonconstant holomorphic function
whose restriction to $A$ is real valued.  If
$f$ vanishes identically on
$A$, then {\rm (3)} holds.  If $f$ does not
vanish identically
 on $A$, but $A$ is a homogeneous CR
submanifold of $\cnn$, then {\rm (3)} still
holds.  
\endproclaim

We shall give another version of conditions
(1) and (2) of Theorem 1, which will give a
reformulation of Theorems 1 and 2. For a  CR
submanifold
$M$ of
$\cnn$,  we  say that  
$M$ is {\it minimal} at $\po \in M$ if there
is no germ of a CR submanifold in $\cnn$
through
$\po$ with the same CR dimension as $M$ at
$\po$, and properly contained in
$M$. A CR submanifold is called {\it generic} at
$p$ if
$$ (T_pM +JT_pM) = T_p \cnn, \tag 0.1
$$ where $T_p\cnn$ is the real tangent space
of $\cnn$.  (See \S 1.1 for more details and
equivalent formulations.)

For an irreducible real algebraic subset $A$
of $\cnn$, we let
$\acr$ be the subset of points in $\areg$ at
which $A$ is CR.  The following contains
Theorems 1 and 2.  
\proclaim {Theorem 3} Let
$A
\subset
\cnn$ be an irreducible real algebraic set,
and let {\rm (1)}, {\rm (2)}, and {\rm (3)}
be the conditions of Theorems {\rm 1} and
{\rm 2}.  Consider also the following
conditions.  
\roster
\item  "(i)" There exists $p \in
\acr$ at which
$A$ is holomorphically nondegenerate. 
\item  "(ii)"  There exists $p \in
\acr$ at which $A$ is generic.
\item "(iii)" There exists $p \in
\acr$ at which $A$ is minimal.
\endroster  Then condition {\rm(i)} is
equivalent to condition {\rm (1)}, and
conditions {\rm(ii)} and {\rm(iii)} together
are equivalent to condition {\rm (2)}. In
particular, {\rm(i)}, {\rm(ii)}, and
{\rm(iii)} together imply the conclusion of
Theorem {\rm 1}.  If either {\rm(i)} or
{\rm(ii)} does not hold, then {\rm (3)} must
hold.  If {\rm(iii)} does not hold, and
$A$ is a homogeneous CR manifold, then {\rm
(3)} must also hold.  
\endproclaim

Note that conditions (i), (ii), and (iii) of
Theorem 3 are all independent of each
other.  

If $M$ is a real analytic CR submanifold of
$\cnn$ and $\po \in M$ (with $M$ not
necessarily minimal at
$\po$), then by Nagano's theorem [N] there
exists a real analytic minimal CR
submanifold of
$M$ through $\po$ of minimum possible
dimension (and the same CR dimension as $M$)
contained in $M$.  Such a manifold is called
the {\it CR orbit} of
$\po$. We call  the germ of the smallest
complex analytic manifold of $\cnn$
containing the CR orbit the {\it intrinsic
complexification} of this orbit. 

Note that if $\scrV \subset \cnn$ is a
complex algebraic set, i.e. defined by the
vanishing of holomorphic polynomials, then
one can define the notion of an algebraic
holomorphic function on an open subset of
$\scrV_{\text {\rm reg}}$ (see \S 3.1).

For point CR submanifolds which are
nowhere minimal, we have the following.

\proclaim {Theorem 4} Let $M$ be a real
algebraic CR submanifold of
$\cnn$ and $\po \in M$. Then the CR orbit of
$\po$ is a real algebraic submanifold of $M$
and  its intrinsic
 complexification, $X$, is a complex
algebraic submanifold of
$\cnn$. For any germ $H$ of a biholomorphism
at
$\po$ of $\cnn$ into itself mapping
$M$ into another real algebraic  manifold of
the same dimension as that of $M$, the
restriction of $H$ to
$X$ is algebraic.   
\endproclaim

The algebraicity of the mapping in Theorem 4
follows from Theorem 1, after it is shown,
in the first part of the theorem, that the
CR orbits are algebraic. (See Theorem
2.2.1.) We mention here that the algebraic
analog of the Frobenius or Nagano theorem
does not hold, since the integral curves of
a vector field with algebraic coefficients
need not be algebraic.  It is therefore
surprising that the CR orbits of an
algebraic CR manifold are algebraic.     In
\S 3.1 we formulate and prove  Theorem
3.1.2, a more general result containing
Theorems 1 and  4, which also applies to points in
an algebraic set
$A$ at which $A$ is not necessarily CR or
even regular, and which, in some cases,
yields algebraicity on a larger submanifold
than the one obtained in Theorem 4.  (See 
Example 3.1.5.)

Note that if a germ of a holomorphic
function is algebraic, it extends as a
(multi-valued) holomorphic function in all of
$\cnn$ outside a proper complex algebraic
subset.  This may be viewed as one of the
motivations for proving algebraicity of
functions and mappings.

We give here a brief history of some
previous work on the algebraicity of
holomorphic mappings between real algebraic
sets. Early in this century Poincar\'e [P]
proved that if a biholomorphism defined in
an open set in $\bC^2$ maps an open piece of
a sphere into another, it
 is necessarily a rational map.  This result
was extended by Tanaka [Ta] to spheres in
higher dimensions. Webster [W1] proved a
far-reaching result for algebraic,
Levi-nondegenerate real hypersurfaces in
$\cnn$;  he proved that any biholomorphism
mapping such a hypersurfaces into another is
algebraic.  Later, Webster's result was
extended in some cases to Levi-nondegenerate
hypersurfaces in complex spaces of different
dimensions (see e.g. Webster [W2],
Forstneri\v c [Fo], Huang [H] and their
references). See also Bedford-Bell [BB] for
other results related to this work. We refer the reader in
addition to the work of Tumanov and Henkin [TH] and
Tumanov [Tu2] which contain results on mappings of higher
codimensional quadratic manifolds. See also related
results of Sharipov and Sukhov [SS] using Levi form
criteria;  some of these results are special cases of
the present work.

  It should be perhaps mentioned  that
the algebraicity results here are deduced from
local analyticity in contrast with the
general ``G.A.G.A. principle" of Serre [Ser], which
deals with the algebraicity of global analytic
objects.

The results and techniques in the papers
mentioned above have been applied to other
questions concerning mappings between
hypersurfaces and manifolds of higher
codimension. We mention here, for instance,
the classification of ellipsoids in $\cnn$
proved in [W1] (see also [W3] for related
problems). We refer also to the regularity
results for CR mappings, proved in Huang
[H], as well as the recent joint work of
Huang with the first and third authors [BHR].
Applications of the results and techniques
of the present paper to CR automorphisms of
real analytic manifolds of higher codimension
and other related questions will be given in
a forthcoming paper of the authors
[BER].     

In [BR3], the first and third authors
proved that for real algebraic hypersurfaces
in
$\cnn$, $N > 1$, holomorphic nondegeneracy
is a necessary and sufficient condition for
algebraicity of all biholomorphisms between
such hypersurfaces.  It should be noted that
any real smooth hypersurface $M
\subset \cnn $ is CR at all its points, and  
if such an $M$ is real analytic and
holomorphically nondegenerate (and $N > 1$), it is
minimal at all points outside a proper analytic
subset of $M$.  Hence, the main result of
[BR3] is contained in Theorem 3 above.  (In
fact the proofs given in this paper are, for
the case of a hypersurface, slightly
simplified from that in [BR3], see [BR4].)
It is easy to check that in
$\bC$, any real algebraic hypersurface 
(i.e. curve) is holomorphically
nondegenerate, but never minimal at any
point.  In fact, by the (algebraic) implicit
function theorem, such a curve is locally
algebraically equivalent to the real line,
which is a homogeneous algebraic set in the
sense of Theorem 3. The conclusion of Theorem
3 agrees with the observation that, for
instance, the mapping $Z\mapsto e^Z$ maps
the real line into itself.

The definition of holomorphic degeneracy was
first introduced by Stanton [St1] for the
case of a hypersurface.  It is proved in
[BR3] (see also [St2]) that if $M$ is a
connected real analytic hypersurface, then
$M$ is holomorphically degenerate at one
point if and only if $M$ is holomorphically
degenerate at all points.  This condition is
also equivalent to the condition that $M$ is
nowhere essentially finite (see
\S 1). In higher codimension we show in this
paper that holomorphic degeneracy propagates
at all CR points (see \S 1.2). The
definition of minimality given here was
first introduced by Tumanov [Tu1].  For real
analytic CR manifolds minimality is
equivalent (by Nagano's theorem [N]) to
the finite type condition of Bloom-Graham
[BG] (see also [BR1]). Both
formulations, i.e. minimality and finite
type, are used in this paper.

 The main technical novelty of this work is
the use of a sequence of sets, called here
the {\it Segre sets} attached to every point
in a real analytic CR manifold. For $M$
algebraic, the Segre sets are (pieces of)
complex algebraic varieties. Another result
of this paper, of independent interest, is a
new characterization of minimality (or
finite type) in terms of Segre sets (see
Theorem 2.2.1).  In fact, it is shown that
the largest Segre set attached to a point
$\po
\in M$ is the intrinsic complexification of
the CR orbit of $\po$. This in particular
proves the algebraicity of the CR orbit when
$M$ is algebraic. The first Segre set of a
point coincides with the so-called Segre
surface introduced by Segre [Seg] and used in
the work of Webster [W1], Diederich-Webster
[DW], Diederich-Fornaess [DF] and others.
Our subsequent Segre sets are all unions of
Segre surfaces. The difficulty in the present
context arises from the fact that the real
algebraic sets considered can be of real
codimension greater than one. Indeed, in the
codimension one case, i.e. hypersurface, the
Segre sets we construct reduce to either the
classical Segre surfaces or to all of
$\cnn$.

The paper is organized as follows.  In \S
1.1 we recall some of the basic definitions
concerned with real analytic manifolds in
$\cnn$ and their CR structures.  The other
subsections of \S 1 are devoted to proving
the main properties of holomorphic
nondegeneracy, which are crucial for the
proofs of the results of this paper.  In \S
2 we introduce the notion of Segre sets, as
described above; their basic properties,
including the characterization of finite
type and the algebraicity of the CR orbits,
are given in Theorem 2.2.1. In \S 3 we prove
the main results of this paper, of which
Theorems 1--4 are consequences.  For the
proof of the most inclusive result, Theorem
3.1.2, a general lemma on propagation of
algebraicity, which may be new, is needed; 
it is proved in \S 3.2. The actual proofs
of Theorems 1--4 are given in \S 3.6.  Examples are
given throughout the paper.

\heading 1. Holomorphic nondegeneracy of real
analytic sets \endheading

\subhead 1.1 Preliminaries on real
submanifolds of $\cnn$ \endsubhead

Let $M$ be a real analytic submanifold of
$\cnn$ of codimension
$d$ and
$\po
\in M$ . Then $M$ near $\po$ is given by
$\r_j(Z, \overline Z) = 0$, $j =1,\ldots,d$, where
the $\r_j$ are real analytic, real-valued
functions satisfying
$$ d\r_1(Z,\overline Z)\wedge...\wedge
d\r_d(Z,\overline Z)\not=0$$ for $Z$ near
$\po$. It can be easily checked that  the
manifold
 $M$ is CR at $\po$ if, in addition, the
rank of $ (\pa\r_1(Z,\overline Z), \ldots,
\pa\r_d(Z,\overline Z))$ is constant for $Z$ near
$\po$, where $\pa f =
\sum_j{\pa f
\over \pa Z_j } dZ_j$.  Also,  $M$ is
generic at $\po$ if the stronger
condition
$$ \pa\r_1(Z,\overline Z)\wedge...\wedge
\pa\r_d(Z,\overline Z)\not=0 \tag 1.1.1$$ holds for
$Z$ near
$\po$. 

For $p \in M$, we denote by
$T_pM$ the real tangent space of $M$ at
$p$ and by $\bC T_pM$ its complexification. 
We denote by
$T^{0,1}_pM$ the complex subspace of
$\bC T_pM$ consisting of all antiholomorphic
vectors tangent to
$M$ at $p$, and by $T^c_p M = \re
T^{0,1}_pM$ the complex tangent space of $M$
at $p$ considered as a real subspace of
$T_pM$.  If $M$ is CR, then dim$_\bC
T^{0,1}_p M$  and dim$_\bR T^c_pM$  are
constant, i.e. independent of
$p$, and we denote by  $T^{0,1}M$ and
$T^c M$ the associated bundles.  The {\it CR
dimension} of $M$ is then
$$\text {\rm CRdim}\ M=\dim_\bC T^{0,1}_p M
=1/2 \, \dim_\bR T^c_pM.$$   If
$M$ is generic, then dim$_\bC T^{0,1}_p M=
N-d$ for all $p$.  If $M$ is CR, then by
Nagano's theorem [N]
$M$ is the disjoint union of real analytic
submanifolds, called the {\it CR orbits} of
$M$. The tangent space of such a submanifold
at every point consists of the restrictions
to that point of the Lie algebra generated
by the sections of
$T^c M$. Hence $M$ is of 
 finite type (in the sense of Bloom-Graham [BG])
or  minimal at
$p$ (as defined in the introduction) if the
codimension of the CR orbit through $p$ is
$0$, i.e. if the Lie algebra generated by
the  sections of
$T^cM$ spans the tangent space of
$M$ at
$p$.

Note that if $M$ is a real analytic 
submanifold 
 of $\cnn$ then there is a proper real
analytic subvariety $V$ of $M$ such that
$M \backslash V$ is a CR manifold. If
$M$ is CR at $\po$ then we may find local 
coordinates
$Z=(Z',Z'')$ such that near $\po$,
$M$ is generic in the subspace
$Z''=0$.  Hence, any real analytic CR
manifold $M$ is a generic manifold in a
complex holomorphic submanifold $\Cal X$ of
$\cnn$, here called the {\it intrinsic
complexification} of $M$. We call
$\dim_\bC \Cal X - \text {\rm CRdim}\ M$ the
{\it CR codimension } of $M$. Hence, if $M$
is a generic submanifold of $\cnn$ of
codimension $d$ its CR dimension is
$N-d$ and its CR codimension is $d$. In view
of the observation above, we shall restrict
most of our analysis to that of generic
submanifolds of
$\cnn$.

For a CR manifold $M$, we define its {\it
H\"ormander numbers} at $\po
\in M$ as follows.  We let $E_0 = T^c_{\po}
M$ and $\mu_1$ the smallest
 integer $\ge 2$ such that the sections of
$T^cM$ and their commutators of lengths $\le
\mu_1$ evaluated at $\po$ span a  subspace
$E_1$ of $T_{\po} M$ strictly bigger than
$E_0$.  The {\it multiplicity} of the first
H\"ormander number
$\mu_1$ is then $\ell_1 = \dim_\bR E_1 - 
\dim_\bR E_0$. Similarly, we define
$\mu_2$ as the smallest integer such that
the sections of of $T^cM$ and their
commutators of lengths $\le
\mu_2$ evaluated at $\po$ span a subspace
$E_2$ of $T_{\po} M$ strictly bigger than
$E_1$, and we let $\ell_2 = 
\dim_\bR E_2 -  \dim_\bR E_1$ be the
multiplicity of
$\mu_2$.  We continue inductively to find
integers $2 \le \mu_1 < \mu_2 <
\ldots < \mu_s$, and subspaces
$T^c_{\po}M = E_0 \subsetneq E_1
\subset
\ldots \subsetneq E_s \subset T_{\po}M$,
where $E_s$ is the subspace spanned by the
Lie algebra of the sections of
$T^c M$ evaluated at $\po$. The multiplicity
$\ell_j$ of each $\mu_j$ is defined in the
obvious way as above.  It is convenient to
denote by $m_1 \le m_2
\le \dots \le m_r $ the {\it H\"ormander
numbers with multiplicity} by taking
$m_1=m_2=\dots=m_{\ell_1}=\mu_1$, and so on.
Note that if $M$ is generic, then
$r=d$ if and only if $M$ is of finite type at
$\po$. More generally, if $M$ is CR, then
$r$ coincides with the CR codimension of $M$
if and only if
$M$ is of finite type at $\po$.

Now suppose that $M$ is a real analytic
generic submanifold of codimension $d$ in
$\cnn$ and
$\rho(Z,\overline Z)=
(\r_1(Z,\overline Z),\ldots,\r_d(Z,\overline Z))$ is a
defining function for
$M$ near $\po \in M $. We write  $N = n+d$.
We define the germ of an analytic subset
$\scrV_{\po}\subset  \cnn$ through
$\po$ by
$$
\scrV_{\po} = \{Z: \r(Z,\z) = 0 \
\hbox {for all} \ \z \ \hbox {near}\
\overline \po\ \hbox {with}\ \r(\po, \z) =0\}.
\tag 1.1.2
$$ Note in fact that $\scrV_{\po}
\subset M$. Then $M$ is called {\it
essentially finite} at
$\po$ if $\scrV_{\po} = \{\po\}$. 

 Recall that by the use of the implicit
function theorem (see [CM], [BJT], [BR2])
we can find holomorphic coordinates $(z,w), z
\in
\cn, w
\in \cd$ vanishing at $\po$ such that near
$\po$, 
$$
\rho(Z,\bar Z)=\im w-\phi(z,\bar z,\re w),
$$ where $\phi(z,\bar z,s)=(\phi_1(z,\bar
z,s),...,\phi_d(z,\bar z,s))$ are 
real-valued real analytic functions in
$\bR^{ 2n+d}$ extending as holomorphic
functions
$\phi(z,\chi,\sigma)$ in
$\bC^{2n+d}$ with 
$$
\phi(z,0,\sigma)\equiv\phi
(0,\chi,\sigma)\equiv0.
$$ Hence, solving in $w$ or $\overline w$ we can
write the equation of $M$ as  
$$ w = Q(z,\overline z,\overline w)\ \ {\text {\rm or}} \ \
\overline w = \overline Q(\overline z,z,w), \tag 1.1.3 
$$ where $Q(z,\ch,\t)$ is holomorphic in a
neighborhood of
$0$ in
$\bC^{2n+d}$, valued in $\cd$ and satisfies
$$Q(z,0,\t) \equiv Q(0,\ch,\t)
\equiv\t . \tag 1.1.4
$$
 It follows from the reality of the
$\rho_j$ and (1.1.3) that the following
identity holds for all $z,
\ch, w \in C^{2n+d}$ near the origin: 
$$ Q(z,\ch,\overline Q(\ch,z,w)) \equiv w.
\tag 1.1.5
$$ Coordinates $(z,w)$ satisfying the above
properties are
 called {\it normal coordinates} at
$\po$.  

If $Z=(z,w)$ are normal coordinates at
$\po$, then the analytic variety defined in
(1.1.2) is given by
$$
\scrV_{\po} = \{ (z,0): Q(z,\ch,0) = 0
\
\hbox {for all } \ \ch \in
\cnn\}. \tag 1.1.6
$$

\remark {Remark {\rm 1.1.1}} If the generic
submanifold $M$ is real algebraic,  then after a
holomorphic algebraic change of coordinates
one can find normal coordinates
$(z,w)$ as above such that the function $Q$
in (1.1.3) is algebraic holomorphic in a
neighborhood of $0$ in
$\bC^{2n+d}$, and hence
$\scrV_{\po}$ is a complex algebraic
manifold. If $M$ is a real algebraic CR 
submanifold, then its intrinsic complexification is a
complex algebraic submanifold. Indeed, these are
obtained by the use of the implicit function theorem,
which preserves algebraicity. (See [BM] and [BR3]
for more details.) 
\endremark  

\subhead 1.2 Holomorphic nondegeneracy and
its propagation \endsubhead

A real analytic  submanifold
$M$ of $\cnn$ is called {\it holomorphically
degenerate} at $\po
\in M$ if there exists a  vector field $X =
\sum_{j=1}^N a_j(Z){\pa \over \pa Z_j}$
tangent to $M$ where the
$a_j(Z)$ are germs of holomorphic functions
at
$\po$ not all vanishing identically on $M$. 
For CR submanifolds,  we shall show that
holomorphic nondegeneracy is in fact
independent of the choice of the point
$\po$.  
\proclaim {Proposition 1.2.1} Let 
$M$ be a connected real analytic CR
submanifold of $\cnn$, and let $p_1, p_2 \in
M$.
 Then
$M$ is holomorphically degenerate at
$p_1$ if and only if it is holomorphically
degenerate at $p_2$.
\endproclaim
\demo {Proof} Since, as observed in
\S 1.1, every CR manifold is a generic
submanifold of a complex manifold, it
suffices to assume that $M$ is a generic
submanifold of $\cnn$. We shall be brief
here, since the proof is very similar to
that of the case where $M$ is a
hypersurface, i.e.
$d=1$, given in [BR3].  We start with an
arbitrary point $\po
\in M$ and we choose normal coordinates
$(z, w)$ vanishing  at
$\po$.  We assume that $M$ is given by
(1.1.3) for $(z, w)$ near
$0$.
 We  write
$$
\overline Q(\ch,z,w) = \sum_\a q_\a(z,w)\ch^\a
\tag 1.2.1
$$  for $ |z|, |\ch|,|w| < \delta$.  We
shall assume that $\delta$ is chosen
sufficiently small so that the right hand side
of (1.2.1) is absolutely convergent. Here
$q_\a$ is a holomorphic function
 defined for $|z|,|w| < \d$ valued in
$\cd$. We leave the proof of the following
claim  to the reader, since it is
very similar to the case
$d=1$ proved in [BR3]:

{\sl Let $(z^1,w^1) \in M $, with $
|z^1|, |w^1| < \delta$.   If
$X$ is a germ at $(z^1,w^1)$ of a
holomorphic vector field in
$\cnn$, then
$X $ is tangent to $M$ if and only if 
$$ X =\sum_{j=1}^na_j(z,w){\pa \over
\pa z_j} \ \ \text {\rm and}   \  \
\sum_{j=1}^na_j(z,w) q_{\a,{z_j}}(z,w)
\equiv 0   , \tag 1.2.2
$$ with $a_j$ holomorphic in a neighborhood
of $(z^1,w^1)$, for all multi-indices
$\alpha$, and
$(z,w) $ in a neighborhood of
$(z^1,w^1)$, where the $q_{\a,z_j}$ are the
derivatives with respect to $z_j$ of the $q_\a$  given
by {\rm (1.2.1)}.
}

As in [BR3], it easily follows by linear
algebra from (1.2.2) that if
$M$ is holomorphically degenerate at a point
$(z^1,w^1)$ as above, then it is
holomorphically degenerate at any point
$(z,w)$ in the local chart of normal
coordinates.  Proposition 1.2.1 then follows
by the existence of normal coordinates at
every point and the connectedness of $M$.
$\square$ 
\enddemo

In view of Proposition 1.2.1, if $M$ is a
connected CR manifold in
$\cnn$ we shall say that $M$ is {\it
holomorphically nondegenerate} if it is
holomorphically nondegenerate at some point,
and hence at every point, of $M$.

\subhead 1.3 The Levi number and essential
finiteness 
\endsubhead

  Let $M$ be a real analytic generic
manifold in $\Bbb C^{N }$,
$\po
\in M$ and
$\rho(Z,\overline Z)$  defining functions for
$M$ near $\po$ as in (1.1.1). Without loss
of generality, we may assume
$\po = 0$. For
$p_1$ close to $0$ we define the 
manifold
$ \Sigma_{p_1 }$ by
$$ \Sigma_{p_1 }= \{ \z \in \Bbb C^{ N}:
\rho(p_1, \z) = 0 \}.
$$ (This is the complex conjugate of the classical Segre
manifold.) Note that by (1.1.1),
$\Sigma_{p_1 }$ is a germ of a smooth
holomorphic manifold in
$\cnn$ of codimension $d$. Let $L_1,
\ldots, L_n$,
$n = N-d$,  given by $L_j =
\sum_{k=1}^Na_{jk}(Z,\overline Z)
\pa / \pa {\overline Z_k} $, be a basis of the CR
vector fields on
$M$ near
$0$ with the $a_{jk}$ real analytic (i.e. a
basis near $0$ of the sections of the bundle
$T^{0,1}M$.)   If
$X_1,
\ldots, X_n$
 are the complex vector fields given by
$X_j =
\sum_{k=1}^Na_{jk}(p_1,\z)
\pa/\pa {\z_k} $, $j =1,\ldots, n$, then
$X_j$ is tangent to
$\Sigma_{p_1 }$ and the $X_j$ span the
tangent space to
$\Sigma_{p_1 }$ for $\z \in
\Sigma_{p_1 }$ in a neighborhood of
$0$, with
$(p_1,\z)
\mapsto a_{jk}(p_1,\z)$ holomorphic near
$(0,0)$ in $\Bbb C^{2N }$.    For a
multi-index
$\alpha = (\a_1,\ldots,\a_n)$ and
$j = 1,\ldots,d$, we define
$c_{j\a}(Z,p_1,\z)$ in $\Bbb C\{Z,p_1,\z\}$,
the ring of convergent power series in $3N$
complex variables, by 
$$ c_{j\a}(Z,p_1,\z)= X^{\a}
\rho_j(Z+p_1, \z), \ \ \ j=1,\ldots,d,
\tag 1.3.1
$$ where $X^\a = X_1^{\a_1}\cdots
X_n^{\a_n}$.

Note that since the $X_j$ are tangent to
$\Sigma_{p_1 }$, we have
$c_{j\a}(0,p_1,\z) = 0$  for all
$(p_1 ,\z )$ near $(0,0)$ and
$\z \in
\Sigma_{p_1 }$. In particular,
$c_{j\a}(0,p_1,\overline p_1) = 0$ for $p_1
\in M$ close to $0$.  It can be checked that
$M$ is  essentially finite at
$p_1$ if the functions $Z \mapsto c_{j\a}(Z,
p_1, \overline p_1)$, $1\le j\le d$,
$\a \in \Bbb Z_+^n$, have only $0$ as a
common zero near the origin for
$p_1$ fixed, small.  (See [BR2] or [BHR]
for a similar argument in the case of a
hypersurface.)

For $1\le j\le d$, $\a \in \Bbb Z_+^n$, let
$V_{j\a}$ be the real analytic $\cnn$-valued functions 
defined near $0$ in
$\Bbb C^{N }$ by 
$$ V_{j\a} (Z,\overline Z)= L^\a
\rho_{jZ}(Z,\overline Z), \tag 1.3.2
$$ where $\rho_{jZ}$ denotes the gradient of
$\rho_j$ with respect to
$Z$ and $L^\a = L_1^{\a_1}\ldots
L_1^{\a_n}$, where
$L_1,\ldots L_n$ are as above.  

In the sequel we shall say that a property
holds {\it generically} on
$M$ if it holds in $M$ outside a proper real analytic
subset. The following definition is independent
of the choice of the defining functions, the
holomorphic coordinates, and the
$L_j$.

 If
$M$ is a generic real analytic submanifold
of $\cnn$ as above, we say that
$M$ is {\it k-nondegenerate} at $Z
\in M$ if the linear span of the vectors
$V_{j\a}(Z,\overline Z)$, $1 \le j \le d, |\a|
\le k$ is all of $\cnn$.  This property is
independent of the choice of the defining
functions $\rho$ and the vector fields $L_j$.

We have the following proposition.
\proclaim {Proposition 1.3.1} Let
$M$ be a connected real analytic generic
manifold of codimension $d$ in $\cnn$. Then
the following conditions are equivalent.
\roster
\item "(i)" $M$ is holomorphically
nondegenerate.
\item "(ii)" There exists $p_1 \in M$ and $k
> 0$ such that $M$ is
$k$-nondegenerate at $p_1$.
\item "(iii)" There exists $V$, a proper
real analytic subset of $M$ and an integer
$\ell= \ell(M) $, $1 \le \ell(M) \le N-d$,
such that $M$ is
$\ell$-nondegenerate at every $p \in
M\backslash V$.
\item "(iv)" There exists $p_1 \in M$ such
that $M$ is essentially finite at $p_1$.
\item "(v)" $M$ is  essentially
finite at all points in a dense open subset of 
$M$.
\endroster
\endproclaim

We shall call the number $\ell(M)$ given in
(iii) above the {\it Levi number} of
$M$.

\demo {Proof} We shall first prove the
equivalence of (i),(ii) and (iii). It is
clear that (iii) implies (ii). We shall now
prove that (ii) implies (i). Assume that
$M$ is
$k$-nondegenerate at $p_1$.  We take normal
coordinates $(z,w)$ vanishing at $p_1$, so
that $M$ is given by (1.1.3)  near
$(z,w) = (0,0)$.  We can take for a basis of
CR vector fields
$$L_j = {\pa \over \pa
\overline z_j}+\sum_{k=1}^d\overline Q_{k\overline z_j}(\overline z,z,w){\pa
\over \pa
\overline w_k}, \ \ j=1,\ldots,n, \tag 1.3.3
$$ so that the $V_{j\a}$ given by (1.3.2)
become, with $Z =(z,w)$,
$$ V_{j\a}(Z,\overline Z)= -\overline Q_{j\overline z^\a
Z}(\overline z,z,w). \tag 1.3.4
$$  The hypothesis (ii) implies that the
vectors $V_{j\a}(0,0)$,
$j=1,\ldots, d$, $|\a| \le k$ span
$\cnn$. By the normality of coordinates,
this implies that the
$q_{j\a z}(0,0)$, $|\a| \le k$, where the
$q_{j\a}(z,w)$ are the components of the
vector
$q_{\a}(z,w)$ defined in (1.2.1), span
$\cn$. This implies, by linear algebra, that
the $a_j(z,w)$ satisfying (1.2.2) in a
neighborhood of $0$ must vanish identically.
Hence $M$ is not holomorphically degenerate
at $0$, proving (i).

Tto show that (i) implies (iii), 
we shall need the following two lemmas, whose
proofs are elementary and left to the
reader.

\proclaim {Lemma 1.3.2} Let
$f_1(\ch),
\ldots,f_d(\ch) $ be $d$ holomorphic
functions defined in an open set
$\Omega$ in
$\Bbb C^{ p}$, valued in
$\Bbb C^{N }$ and generically linearly
independent in $\Omega$.  If the 
$\pa ^\a f_j(\ch)$, $j =1,\ldots, d$,
$\a
\in
\Bbb Z^{p }_+$ span
$\Bbb C^{N }$ generically in
$\Omega$, then the $\partial
 ^\a f_j(\ch)$, $j =1,\ldots, d$,
$|\a| \le N-d$ also span
$\Bbb C^{N }$ generically in
$\Omega$.   
\endproclaim
\proclaim {Lemma 1.3.3} Let $(z,w)$ be
normal coordinates for $M$ as above, and let
$h(\ch,z,w)$ be a holomorphic function in
$2n+d$ variables defined in a connected
neighborhood in
$\bC^{2n+d}$ of
$z=z^1, w=w^1,\ch=\overline z^1$, with
$(z^1,w^1) \in M$, and assume  that
$h(\overline z,z,w)
\equiv 0$, for $(z,w) \in M$.  Then
$h
\equiv 0$.
\endproclaim 

To prove (i) implies (iii), we again take
$(z,w)$ to be normal coordinates around some
point $\po
\in M$.  By the assumption (i) and (1.2.2),
it follows that the $q_{j\a ,z}(z,w)$,
$j=1,\ldots,d$, all $\a$, span $\cn$ 
 generically.  Equivalently, by the
normality of the coordinates, we obtain that
the $ \overline Q_{j\overline z^\a Z}(0,z,w)$ generically
span $\cnn$.  We
claim that the  $\overline Q_{j\overline z^\a Z}(\overline z,z,w)$
generically span $\cnn$ for $(z,w) \in M$.
Indeed, if the
$\overline Q_{j\overline z^\a Z}(\overline z,z,w)$ do not span, then
all
$N\times N$  determinants $\Delta (\overline z,z,w)$
extracted from the components of these
vectors vanish identically on
$M$ and hence, by Lemma 1.3.3,
$\Delta (\ch,z,w)\equiv 0$ in
$\bC^{2n+d}$. In particular, $\Delta
(0,z,w)\equiv 0$, which would contradict 
the fact that the    
$ \overline Q_{j\overline z^\a Z}(0,z,w)$ generically span
$\cnn$. This proves the claim.

Now choose $(z^0,w^0) \in M$ so that
$\Delta(0,z^0,w^0)\not= 0$ for some
determinant $\Delta$ as above.  We apply
Lemma 1.3.2 with
$f_j(\ch)=\overline Q_{jZ}(\ch,z^0,w^0), 
j=1,\dots,d$, to conclude that there exists
$\ell \le N-d$ such that   in the local chart
$(z,w)$, the $V_{j\a}(Z,\overline Z)$ (see (1.3.4))
for
$|\a|\le
\ell$ span $\cnn$ generically for
$Z \in M$.  Since this property is
independent of the choice of local
coordinates, condition (iii) follows from the
connectedness of
$M$. This completes the proof of the
equivalence of (i), (ii) and (iii). 

It remains to show that (i), (ii) and (iii)
are equivalent to (iv) and (v).  We show
first that (iii) implies (iv).  Let
$p_1 \in M$ be any
$\ell$-nondegenerate point, i.e., the span of
$V_{j\a}(p_1,\overline p_1)$, $1 \le j \le d$, $|\a|
\le \ell$, is $\cnn$.  On the other hand, it
follows from (1.3.1) and (1.3.2) that 
$$c_{j\a
Z}(0,p_1,\overline p_1)= 
V_{j\a}(p_1,\overline p_1). \tag 1.3.5$$ 
Hence by the inverse mapping theorem the only
common zero, near
$0$, of the functions $Z \mapsto c_{j\a
}(Z,p_1,\overline p_1)$ is $0$, which proves that
$M$ is essentially finite at $p_1$, hence (iv).

Next, assume that (v) holds. If the
rank of the $V_{j\a}(Z,\overline Z)$ were less than $N$
generically on $M$, then at any point $p_1$ of maximal
rank  near $0$ in $M$, in view of
(1.3.5) and the implicit function theorem, there
would exist a complex curve
$Z(t)$ through $0$ such that $c_{j\a}(Z(t),p_1,\overline p_1)=0$
for all small $t$ and all $j, \a$.  Hence $M$ would not
be essentially finite at $p_1$, contradicting (v),
since $p_1$ can be chosen in an open dense set.

Since (v) implies (iv) is trivial, it remains only
to show (iv) implies (v).  For this we need the
following lemma.  

\proclaim {Lemma 1.3.4}  Let $\{f_j\}_{j \in J}$
 be holomorphic in a neighborhood of $0$ in
$\cnn$. Suppose that $Z=0$ is an isolated zero
of the functions $f_j(Z)-f_j(0),\ j\in J$.  Then
there exists
$\d > 0$ such that for $|Z_0| <
\d$,
$Z=0$ is an isolated zero of the functions
$f_j(Z+Z_0) - f_j(Z_0)$, $j \in J$.      
\endproclaim

\demo {Proof}  For $j \in J$, let $F_j(Z,\z) =
f_j(Z)-f_j(\z)$, which is holomorphic near $0$ in
$\bC^{2N}$.  Let $V$ be the variety of zeros of
the $F_j$.  We claim that there exists $\epsilon
> 0$ and $\d > 0$ such that if $|\z_0 | < \d$,
then the set $V\cap \{(Z,\z) \in \bC^{2N}:
|Z|<\epsilon, \z=\z_0\}$ is discrete.  Indeed, by
assumption there exists $\epsilon > 0$ such that 
$V\cap\{|Z|=\epsilon, \z=0\}=\emptyset$. 
Therefore by compactness, there exists $\d$, $0 <\d <
\epsilon$, such that $V \cap \{|Z|=\epsilon, |\z| < \d\}
=
\emptyset$.  Hence for any $|\z_0| < \d$, the set
$V\cap \{|Z| < \epsilon, \z=\z_0\}$ is discrete. 
Hence the zero $Z = \z_0$ of $F(Z,
\z_0)$ is isolated, which completes the proof of
the lemma.  
 $\square$ 
\enddemo     

We may now prove that (iv) implies (v). 
Choose normal coordinates $Z=(z,w)$ around
$p_1 \in M$ at which $M$ is essentially finite, and
observe that if
$\po = (z^0,w^0)$ is in this local chart, we have
$$ c_{j,\a}(Z, \po,\overline p_0) = -\overline Q_{j,\chi^\a}(\overline z^0,
z^0,w^0) + \overline Q_{j,\chi^\a}(\overline z^0,
z+z^0,w+w^0). \tag 1.3.6
$$
By Lemma 1.3.4, we conclude from (1.3.6) that
$M$ is essentially finite for any $\po$ in a
neighborhood of $p_1$.  Property (v) follows by
connectedness of $M$. This completes the proof
of Proposition 1.3.1.   

$\square$ 
\enddemo

\subhead 1.4 Holomorphic
nondegeneracy of real algebraic sets
\endsubhead

Recall that if $A$ is a real algebraic
subset of $\cnn$, we denote by $\areg$ the
set of points at which $A$ is a real
analytic manifold and by $\acr$ the set of
points of
$\areg$ at which $\areg$ is CR. In this
subsection we prove the following result.  

\proclaim {Proposition 1.4.1} Let
$A \subset \cnn$ be an irreducible real
algebraic set and $p_1, p_2
\in A_{CR}$. Then $A$ is holomorphically
degenerate at
$p_1$ if and only if it is holomorphically
degenerate at
$p_2$.
\endproclaim

We note that if $A_{CR}$ is connected then
the proposition follows immediately from
Proposition 1.2.1.  However, even if
$A$ is irreducible, $A$, $\areg$, and $\acr$
need not be connected.

\demo {Proof} It follows from the proof of
Proposition 1.2.1 that if
$M$ is a real algebraic CR manifold,
holomorphically degenerate at $\po \in M$,
then we can find a holomorphic vector field 
$$X = \sum_{j=1}^N a_j(Z){\pa
\over \pa Z_j} \tag 1.4.1
$$ tangent to $M$ with $a_j(Z)$ algebraic
holomorphic near $\po$ and not all vanishing
identically on $M$. Indeed, by Remark 1.1.1,
we may assume that the functions
$\overline Q$ and $q_\a$ in (1.2.1) are algebraic. 
Since the $a_j(z,w)$ in (1.2.2) are obtained
by solving a linear system of equations, we
can find a set of solutions  which are
algebraic.

Assume that $\acr$ is holomorphically
degenerate at
$p_1$.  By the observation above, we can
find $X$ of the form (1.4.1), with the
$a_j(Z)$ holomorphic algebraic, tangent to
$A$ near $p_1$. Since the
$a_j(Z)$ are algebraic, they extend as
multi-valued holomorphic functions to
$\cnn\backslash V$, where $V$ is a proper
complex algebraic subvariety of $\cnn$ with
$p_1 \not\in V$.  Hence $A \cap V$ is a
proper real algebraic subvariety of $A$.  Let
$U$ be a connected open neighborhood of
$p_2$ in $\acr$ and let $p_3 \in U \backslash
V$.  (If $p_2 \not\in V$, we may take $p_3 =
p_2$.) If $d =
\text{\rm codim}_\bR A$, then by a classical
theorem in real algebraic geometry
[HP, Chapter 10], there exist real valued
polynomials
$\rho_1(Z,\overline Z),
\ldots,\rho_d(Z,\overline Z)$ with $A = \{ Z \in
\cnn: \rho_j(Z,\overline Z) = 0, j=1,\ldots,d\}$ and
$d\r_1,\ldots, d\r_d$ generically linearly
independent on $A$. Let $\scrA$ be the
complexification of $A$, i.e. the
irreducible complex algebraic set in
$\bC^{2N}$ given by $\scrA =
\{(Z,\z) \in \ctnn: \rho_j(Z,\z) = 0,
j=1,\ldots,d\}$, and let
$\hat V = V\times \cnn_\z$.  We identify
$\cnn$ with a subset of
$\ctnn$ by the diagonal mapping
$Z \mapsto (Z,\overline Z)$, so that $A$ and $V$
become subsets of $\scrA$ and $\hat V$,
respectively.  We claim that $p_1$ and $p_3$
(considered now as points
$\scrA$) can be connected by a curve
contained in
$\scrA_{\text {\rm reg}}\backslash
\hat V$.  The claim follows from the fact
that 
$\scrA_{\text {\rm reg}}\cap
\hat V$ is a proper algebraic subvariety of
$\scrA$ and hence its complement in $\scrA$
is connected, by the irreducibility of
$\scrA$.  We conclude that the holomorphic
continuation of the vector field (1.4.1), 
thought of as
a vector field in
$\bC^{2N}$, is
tangent to
$\scrA$ at every point along this curve,
from which we conclude that
$A$ is holomorphically degenerate at $p_3$. 
We may now apply Proposition 1.2.1 to the CR
manifold $U$ to conclude that $A$ is also
holomorphically degenerate at $p_2$.  
$\square$ 
\enddemo

\remark  {Remark {\rm 1.4.2}} For a general
real algebraic submanifold
$M \subset
\cnn$, not necessarily CR, it can happen
that $M$ is holomorphically degenerate
at all CR points, but not holomorphically degenerate
at points where $M$ is not CR, as is
illustrated by the following example.
 Let $M
\subset
\bC^4$ be the manifold of dimension
$5$ given by 
$$Z_3 = \overline Z_1^2, \ \ \ \re Z_4 = Z_1\overline Z_2
+Z_2\overline Z_1.$$
$M$ is a CR manifold away from
$Z_1=Z_3=0$, and $M$ and $\mcr$ are
connected. At the CR points the holomorphic
vector fields tangent to
$M$ are all holomorphic multiples of the
vector field $X = \pa /\pa Z_2 + 2
Z_3^{1/2}\pa/\pa Z_4$.  (Note that here
$Z_3^{1/2}=\overline Z_1$ on
$M$.) We conclude that there is no
nontrivial germ of a holomorphic vector
field tangent to $M$ at a non CR  point of
$M$. 
\endremark

\heading 2. The Segre sets
of a real analytic CR  submanifold
\endheading
\subhead 2.1. Complexification of
$M$, involution, and projections\endsubhead
Let
$M$ denote a generic  real analytic
submanifold in some neighborhood
$U\subset \bC^N$ of $p_0\in M$. Let $\r =
(\r_1,\ldots \r_d) $ be defining functions
satisfying (1.1.1) and choose holomorphic
coordinates $Z=(Z_1,\ldots,Z_N)$ vanishing at
$\po$.  Embed 
$\bC^N$ in
$\bC^{2N}=\bC^N_Z\times\bC^N_\zeta$ as the
real plane
$\{(Z,\zeta)\in\bC^{2N}\:\zeta=\bar Z\}$. Let
us denote by
$\text{pr}_Z$ and $\text{pr}_\zeta$ the
projections of $\bC^{2N}$ onto 
$\bC^N_Z$ and $\bC^N_\zeta$, respectively.
The natural anti-holomorphic involution
$\sharp$ in $\bC^{2N}$ defined by
$$ ^\sharp(Z,\zeta)=(\bar\zeta,\bar
Z)\tag2.1.1
$$ 
 leaves the plane
$\{(Z,\zeta)\:\zeta=\bar Z\}$ invariant.
This  involution induces the usual
anti-holomorphic involution in $\bC^N$ by
$$
\bC^N\ni
Z\to\text{pr}_\zeta(^\sharp\text{pr}_Z^{-1}(Z))=\bar
Z\in\bC^N.
\tag2.1.2
$$ Given a set $S$ in $\bC^N_Z$ we denote by
$^* S$ the set in
$\bC^N_\zeta$ defined by
$$
^*S=\text{pr}_\zeta(^\sharp\text{pr}_Z^{-1}(S))=\{\zeta\:\bar\zeta\in
S\}.
\tag2.1.3
$$ By a slight abuse of notation, we use the
same notation for the corresponding 
transformation taking sets in $\bC^N_\zeta$ to
sets in $\bC^N_Z$. Note that if $X$ is a
complex analytic set defined  near
$Z^0$ in some domain $\Omega
\subset\bC_Z^{N}$  by
$h_1(Z)=...=h_k(Z)=0$,  then $^* X$ is the
complex   analytic set in
$^*\Omega\subset\bC^N_\zeta$ defined near
$\zeta^0=\bar Z^0$ by
$\bar h_1(\zeta)=...=\bar h_k(\zeta)=0$.
Here, given a holomorphic function $h(Z)$ we
use the notation 
$
\bar h(Z)=\overline{h(\bar Z)}$. The transformation
$*$ also preserves algebraicity of sets.

Denote by
$\scrM\subset\bC^{2N}$ the complexification
of $M$ given by
$$
\scrM=
\{(Z,\zeta)
\in\bC^{2N}
\:\rho(Z,\zeta)=0\}.\tag2.1.4$$  This is a complex
submanifold of codimension $d$ in some neighborhood
of 
$0$ in $\bC^{2N}$. We choose our neighborhood
$U$ in $\bC^N$ so small that
$U\times{}^*U\subset\bC^{2N}$ is contained in
the neighborhood where
$\scrM$ is a manifold. Note that
$\scrM$ is invariant under the involution 
$\sharp$ defined in (2.1.1). Indeed all the
defining functions
$\rho(Z,\bar Z)$ for $M$ are real-valued,
which implies that the holomorphic extensions
$\rho(Z,\zeta)$ satisfy
$$
\bar\rho(Z,\zeta)=\rho(\zeta,Z).
\tag2.1.5
$$ Thus, given
$(Z,\zeta)\in\bC^{2N}$ we have
$\overline{\rho(^\sharp(Z,\zeta))}=\overline{\rho(\bar\zeta,\bar
Z)}=\bar\rho (\zeta,Z)=\rho(Z,\zeta)$,
 so $^\sharp(Z,\zeta)\in\scrM$ if and only if
$(Z,\zeta)\in \scrM$.

\subhead 2.2. Definition of the Segre sets of $M$
at 
$p_0$\endsubhead We associate to
$M$ at $p_0$ a sequence  of germs of sets
$N_0,N_1,...,N_{j_0}$ at
$p_0$ in $\bC^N$---henceforth  called {\it
the Segre sets} of $M$ at $p_0$ for reasons
that will become apparent---defined as
follows.  Define $N_0=\{p_0\}$ and define the
consecutive sets inductively (the number
$j_0$ will be defined later) by
$$ N_{j+1}=\text{pr}_Z\left(\scrM\cap
\text{pr}_\zeta^{-1}\left(^*N_j\right)
\right)=\text{pr}_Z\left(\scrM\cap{}^\sharp\text{pr}_Z^{-1}\left(N_j\right)
\right).\tag2.2.1
$$ Here, and in what follows, we abuse the
notation slightly by identifying a  germ
$N_j$ with some representative of it. These
sets are, by definition, invariantly defined
and they  arise naturally in the study
of mappings between submanifolds (see \S 3). 

Let the defining functions $\r$ and the
holomorphic coordinates
$Z$ be as in \S 1.1.  Then the sets
$N_j$ can be described as follows, as is
easily verified. For odd
$j=2k+1$ ($k=0,1,...$), we have
$$
\aligned N_{2k+1}=\{Z\: &\exists
Z^1,...,Z^k,\zeta^1,...,\zeta^k\:\\
&\rho(Z,\zeta^k)=\rho(Z^k,\zeta^{k-1})=...=\rho(Z^1,0)=0\,,\\
&\rho(Z^k,\zeta^k)=\rho(Z^{k-1},\zeta^{k-1})=...=\rho(Z^1,\zeta^1)=0\};
\endaligned\tag2.2.2
$$ note that for $k=0$ we  have
$$ N_1=\{Z\:\rho(Z,0)=0\}.\tag2.2.3
$$ For even $j=2k$ ($k=1,2,...$), we have
$$
\aligned N_{2k}=\{Z\: &\exists
Z^1,...,Z^{k-1},\zeta^1,...,\zeta^k\:\\
&\rho(Z,\zeta^k)=\rho(Z^{k-1}
,\zeta^{k-1})=...=\rho(Z^1,
\zeta^1)=0\,,
\\ &\rho(Z^{k-1},\zeta^k)=\rho(Z^{k-2}
,\zeta^{k-1})=...=\rho(0,\zeta^1)= 0\}.
\endaligned\tag2.2.4
$$
 For $k=1$, we have
$$ N_2=\{Z\:\exists
\zeta^1\:\rho(Z,\zeta^1)=0\,,\,\rho(0,\zeta^1)
=0\}.\tag2.2.5
$$
 {}From \thetag{2.2.2} and
\thetag{2.2.4} it is easy to deduce the
inclusions
$$ N_0\subset N_1\subset...\subset
N_j\subset...\tag2.2.6
$$  When $d=1$ the set $N_1$ is the so-called
Segre surface through 0 as introduced by
Segre [S], and used by Webster [W1],
Diederich--Webster [DW],  
Diederich--Fornaess [DF], Chern--Ji [CJ], and
others. Here the set $N_2$ is the union of Segre
manifolds through points $\zeta_1$ such that
$\bar\zeta_1$ belongs to the Segre surface
through 0. Subsequent
$N_j$'s can be described similarily as unions
of Segre manifolds. 

In order to simplify the calculations, it is
convenient to use normal coordinates
$Z=(z,w)$ for
$M$ as in \S1.1. Recall that
$M$ is assumed to be generic and of
codimension $d$; we write $N=n+d$. If $M$ is
given by (1.1.3), it will be convenient to
write
$$ Q(z,\chi,\tau)=\tau+
q(z,\chi,\tau),\tag2.2.7
$$ where 
$$ q(z,0,\tau)\equiv
q(0,\chi,\tau)\equiv0.\tag2.2.8
$$ In $\bC^{2N}$, we choose coordinates
$(Z,\zeta)$ with
$Z=(z,w)$ and
$\zeta=(\chi,\tau)$, where
$z,\chi\in\bC^n$ and
$w,\tau\in\bC^d$. Thus, in view of
\thetag{1.1.3}, the complex manifold 
$\scrM$ is defined by either of the equations
$$ w=Q(z,\chi,\tau)\quad\text{\rm
or}\quad\tau=\bar Q(\chi,z,w).
\tag2.2.9
$$ In normal coordinates, we find that in the
expression
\thetag{2.2.2} for $N_{2k+1}$ we can solve
recursively for
$w^1,\tau^1,w^2,\tau^2,...,w^k,
\tau^k$ and parametrize $N_{2k+1}$ by 
$$
\bC^{(2k+1)n}\ni(z,z^1,...,z^k,\chi^1,...,\chi^k)=\Lambda\mapsto
(z,v^{2k+1}(\Lambda))
\in\bC^N,\tag2.2.10
$$ where 
$$ v^{2k+1}(\Lambda)=
\tau^k+q(z,\chi^k,\tau^k),\tag2.2.11
$$ and recursively
$$
\tau^l=w^l+\bar q(\chi^l,z^l,w^l) \ \ \text
{\rm with}\ \ w^l=\left\{ \aligned &\tau^{l-1}+q(z^{l},
\chi^{l-1},
\tau^{l-1}), \, l\ge 2\\ &0,\,\ l=1
\endaligned\right. \tag2.2.12
$$  for $l=1,2,...,k$; for
$k=0$, we have $v^1\equiv0$. Similarily,  we
can parametrize
$N_{2k}$ by
$$
\bC^{2kn}\ni(z,z^1,...,z^{k-1},\chi^1,...,\chi^k)=\Lambda\mapsto(z,v^{2k}
(\Lambda))
\in\bC^N,\tag2.2.13
$$ where 
$$ v^{2k}(\Lambda)=
\tau^k+q(z,\chi^k,\tau^k),\tag2.2.14
$$ and recursively
$$
\tau^{l+1}=w^{l}+\bar q(\chi
^{l+1},z^{l},w^{l})\ \ \text {\rm with}\ \ 
w^l=\tau^l+q(z^l,\chi^l,\tau^l),
\tag2.2.15
$$ for $l=1,...,k-1$ and
$\tau^1=0$.  Define $d_j$ to be the maximal
rank of the mapping
\thetag{2.2.10} or 
\thetag{2.2.13} (depending on whether $j$ is
odd or even) near
$0\in\bC^{jn}$. It is easy to see that
$d_0=0$ and $d_1=n$. In view of
\thetag{2.2.6}, we  have
$ d_0<d_1\leq d_2\leq d_3\leq\ldots $. We
define the number
$j_0 \ge 1$ to be the greatest integer  such that
we have strict inequalities
$$ d_0<d_1<...<d_{j_0}.\tag2.2.16
$$ Clearly, $j_0$ is a well defined finite
number because, for all $j$, we have
$d_j\leq N=n+d$ and
$d_{j_0}\geq n+j_0-1$ so that we have
$j_0\leq d+1$. The 
$d_j$'s stabilize for $j\geq j_0$, i.e. 
$ d_{j_0}=d_{j_0+1}=d_{j_0+2}=....,
$ by the definition of the Segre sets.

So far we have only considered generic 
submanifolds. If $M$ is a real analytic CR
submanifold of $\bC^N$, then $M$ is generic
as a submanifold of its intrinsic
complexification $\Cal X$ (see \S1.1). If
$M$ is real algebraic then $\Cal X$ is complex
algebraic. The Segre sets of $M$ at a point
$p_0\in M$ can be defined  as subsets of
$\bC^N$  by the process described at the
beginning of this subsection (i.e. by
\thetag{2.2.1}) just as for generic
submanifolds or they can be defined as
subsets of $\Cal X$ by identifying $\Cal X$
near $p_0$ with $\bC^K$ and
considering
$M$ as a generic submanifold of $\bC^K$. It
is an easy exercise (left to the reader) to
show that these definitions are equivalent
(i.e. the latter sets are equal to the former
when viewed as subsets of $\bC^N$).

The main result in this section is the
following. Let the H\"ormander numbers, with
multiplicity, be defined as in \S1.1.
\proclaim{Theorem 2.2.1} Let $M$ be a  real
analytic CR submanifold in
$\cnn$ of CR dimension $n$ and of CR
codimension
$d$ and $p_0\in M$. Assume that there are 
$r$ (finite) H\"ormander numbers of $M$ at
$p_0$, counted with multiplicity. Then the
following hold.
\roster
\item"(a)" There is a holomorphic manifold
$X$ of (complex) dimension
$n+r$ through $p_0$ containing the maximal
Segre set $N_{j_0}$ of
$M$ at $p_0$ (or, more  precisely, every
sufficiently small representative of it) such
that
$N_{j_0}$ contains a relatively open subset
of $X$. In particular, the generic dimension
$d_{j_0}$ of
$N_{j_0}$ equals $n+r$.
\item"(b)" The intersection $M\cap X$ is the
CR orbit of the point
$p_0$ in $M$.
\item"(c)" If $M$ is real algebraic then $X$
is complex algebraic, i.e.
$X$ extends as an irreducible algebraic
variety in
$\bC^N$.\endroster
\endproclaim  In particular, this theorem
gives a new criterion for
$M$ to be of finite type (or minimal) at
$\po$.  The following is an immediate
consequence of the theorem.
\proclaim{Corollary 2.2.2} Let $M$ be a  real
analytic CR submanifold in
$\cnn$ of CR dimension $n$ and of CR
codimension
$d$ and $p_0\in M$. Then
$M$ is minimal at $p_0$,
 if and only if the generic dimension
$d_{j_0}$ of the maximal Segre set
$N_{j_0}$ of $M$ at 
$p_0$ is $n+d$. In particular, if $M$ is
generic, then $M$ is minimal at $\po$ if and
only if $d_{j_0} = N$.
\endproclaim

\noindent{\bf Example 2.2.3.} Let
$M\subset\bC^3$ be the generic submanifold
defined by
$$
\im w_1  =|z|^2, \ \ \  \im w_2
=|z|^4.
$$ 
Then $M$ is of finite type at 0 with
H\"ormander numbers  $2,4$.
 The  Segre sets
$N_1$ and $N_2$ at $0$ are given by
$$ N_1=\{(z,w_1,w_2)\:w_1=0\,,\,w_2=0\}, \tag2.2.17
$$ 
$$
N_2=\{(z,w_1,w_2)\:w_1=2iz\chi\,
,\,w_2=2iz^2\chi^2\,,\,\chi\in\bC\}.\tag2.2.18
$$ Solving for $\chi$ in (2.1.18) we obtain in
this way (outside the plane $\{z=0\}$)
$$
N_2=\{(z,w_1,w_2)\:w_2=-iw_1^2/2\}.
$$ Using the definition
\thetag{2.2.1}, we obtain
$$ N_3=\{(z,w_1,w_2)\:
w_2=iw_1(w_1/2-2z\chi)\,,\,\chi\in\bC\}.
$$
 We have $d_3=3$; $N_3$ contains $\bC^3$
minus the planes
$\{z=0\}$ and $\{w_1=0\}$.
\medskip
\noindent{\bf Example 2.2.4.} Consider
$M\subset\bC^3$ defined by
$$ \im w_1 =|z|^2,\ \ \ \im w_2=\re w_2|z|^4.
$$ Here $2$ is the only H\"ormander
number at the origin.  Again, $N_1$ is given
by
\thetag{2.2.17}, and
$$
N_2=\{(z,w_1,w_2)\: z \not=
0,w_2=0\}\cup \{0,0,0\}.
$$ It is easy to see that subsequent Segre
sets are equal to
$N_2$. Thus, $N_2$ is the maximal Segre set
of $M$ at 0 , $d_2=2$, and the intersection of
(the closure of) $N_2$ with $M$ equals the
CR orbit  of $0$.\medskip

Let us also note that part (c) of Theorem
2.2.1 implies the following.
\proclaim{Corollary 2.2.5} The CR orbits of a
real algebraic CR manifold are
algebraic.\endproclaim The theorem of Nagano
([N]) states that the integral manifolds of
 systems of vector fields, with real
analytic coefficients, are real
analytic. Thus,  the CR orbits of a real
analytic CR manifold $M$
are real analytic submanifolds of $M$.
However, in
general the integral manifolds of systems of
vector fields with  real algebraic
coefficients are not algebraic manifolds, as
can be readily seen by examples. Hence, one
cannot use Nagano's theorem to deduce that the
orbits of an algebraic CR manifold are
algebraic. Corollary 2.2.5 seems not to have
been known before. 

Before we prove Theorem 2.2.1 (in
\S2.5) we first discuss the homogeneous case
because the proof of the theorem will
essentially reduce to this case. We first
consider the case where the CR dimension is 1
(\S2.3) and then give the modifications
needed to consider  the general case (\S2.4).

\subhead 2.3. Homogeneous submanifolds of CR
dimension 1\endsubhead
Let $\mu_1 \le \ldots \le\mu_N$ be $N$
positive integers. For $t > 0$ and 
$Z=(Z_1,\ldots,Z_N)\in\cnn$, we let $\d_tZ
= (t^{\mu_1}Z_1, \ldots,t^{\mu_N}Z_N)$. A
polynomial $P(Z,\overline Z)$ is {\it weighted
homogeneous of degree $m$} with respect to the
weights   
$\mu_1, \ldots ,\mu_N$ if $P(\d_tZ,\d_t\overline Z)
= t^mP(Z,\overline Z)$ for $t > 0$.
 
In this  section and
the next, we consider submanifolds
$M$ in $\bC^{N}$, $N = n+d$, of the form
$$ M\:\left\{\aligned w_1 &=\bar
w_1+q_1(z,\bar z)\\&\vdots\\w_j &=\bar w_j+
q_j(z,\bar z,\bar w_1,...,\bar
w_{j-1})\\&\vdots\\w_r &=\bar w_r+ q_r(z,\bar
z,\bar w_1,...,\bar w_{r-1})\\w_{r+1} &=\bar
w_{r+1}\\&\vdots\\ w_d &=\bar
w_d,\endaligned\right.\tag2.3.1
$$ where $0\leq r\leq d$ is an integer ($r=0$
corresponds to the canonically flat
submanifold), and each $q_j$, for
$j=1,...,r$, is a weighted homogeneous
polynomial of degree
 $m_j$. The weight of each $z_j$ is 1 and the
weight of $w_k$, for
$k=1,...,r$, is $m_k$. Since  the defining
equations of $M$ are polynomials, we can, and
we will, consider the sets
$N_0,...,N_{j_0}$ attached to $M$ at 0 as
globally defined subsets of $\cnn$. Each 
$N_j$ is contained in an irreducible complex
algebraic variety of dimension $d_j$ (here,
an algebraic variety of dimension $N$ is the
whole space $\bC^{N}$). The latter follows
from the parametric definitions
\thetag{2.2.10} and
\thetag{2.2.13} of $N_j$ and the algebraic
implicit function theorem.

We let $\pi_j$, for $j=2,..., d+1$, be the 
projection
$\pi_j\:\bC^{n+d}\mapsto\bC^{n+j-1}$ defined
by 
$$
\pi_j(z,w_1,...w_d)=(z,w_1,...,w_{j-1})
.\tag2.3.2
$$ We define
$M^j\subset\bC^{n+j-1}$ to be
$\pi_j(M)$. By the form 
\thetag{2.3.1} of
$M$, it follows that each $M^j$ is the CR
manifold of codimension
$j-1$  defined by the $j-1$ first equations
of \thetag{2.3.1}. Throughout this  section
and the next, we work under the assumption
that $M$ satifies the  following.
\proclaim{Condition 2.3.1} The CR manifold
$M^j$, for $j=2,...,r+1$,  is of finite type
at 0.\endproclaim 

For clarity, we consider first the case
where the CR dimension, $n$, is one, i.e.
$z\in\bC$.  The rest of this section is
devoted to this case.
\proclaim{Proposition 2.3.2} Let
$M$ be of the form \thetag{2.3.1} with CR
dimension $n=1$ and assume that $M$ satisfies
Condition {\rm 2.3.1}. Let 
$N_0,N_1,...,N_{j_0}$ be the Segre sets
 of $M$ at $0$, and let
$d_0,d_1,...,d_{j_0}$ be their generic
dimensions. Then $j_0=r+1$ and $d_j=j$, for
$0\leq j\leq r+1$. Furthermore, for each
$j=0,...,r+1$, there is a proper
complex algebraic variety
$V_j\subset\bC^{j}$ such that $N_j$ satisfies
$$
\aligned&N_j\cap\left((\bC^j\setminus
V_j)\times\bC^{d-j+1}\right)=\\
&\left\{(z,w_1,...,w_d)\in\left((\bC^j
\setminus
V_j)\times\bC^{d-j+1}\right)
\: 
w_k=f_{jk}(z,w_1,...,w_{j-1})\,,\,
k=j,...,d\right\},\endaligned\tag2.3.3
$$ where each $f_{jk}$, for
$k=j,...,r$,  is a (multi-valued) algebraic
function  with $b_{jk}$ holomorphic, disjoint
branches outside $V_j$ and where $f_{jk}
\equiv0$ for
$k=r+1,...,d$.\endproclaim
\demo{Proof} Clearly, the first statement of
the proposition follows from  the last one.
Thus, it suffices to prove that, for each
$j=0,...,r+1$, there is a proper algebraic
variety $V_j$ such that
\thetag{2.3.3} holds. The proof of this is by
induction on $j$. 

Since $N_0=\{0\}$ and
$N_1=\{(z,w)\:w=0\}$,
\thetag{2.3.3} holds for $j=0,1$ with
$V_0=V_1=\emptyset$. We
assume that there are
$V_0,...,V_{l-1}$ such that
\thetag{2.3.3} holds for
$j=0,...,l-1$. By
\thetag{2.2.1}, we have
$$
N_l=\left\{(z,w)\:\exists(\chi,\tau)\in{}^*N_{l-1}\,,\,
(z,w,\chi,
\tau)\in\scrM
\right\}.\tag2.3.4
$$
\proclaim{Assertion 2.3.3} The set of points
$(z,w_1,...,w_{l-1})\in\bC^{l}$  such that
there exists
$(w_l,...,w_d)\in\bC^{d-l+1}$ and 
$(\chi,\tau)\in{}^*(N_{l-1}\cap(V_{l-1}\times\bC^{d-l+2}))$
with the property that
$(z,w,\chi,\tau)\in\scrM$  is contained in a
proper algebraic variety
$A_l\subset\bC^{l}$.\endproclaim
\demo{Proof of Assertion 2.3.3} Let $S$ be the 
set of points
$(z,w_1,..., w_{l-1})\in\bC^l$ described in
the assertion. Then
$(z,w_1,...,w_{l-1})\in\bC^l$ is in $S$ if 
$$
\tau_j=w_j+\bar
q_j(\chi,z,w_1,...,w_{j-1})\quad,\quad
j=1,...,l-1.
\tag2.3.5
$$ 
for some
$(\chi,\tau_1,...,\tau_{l-1})\in{}^*(\pi_l
(N_{l-1})\cap(V_{l-1}\times\bC))$. (Recall the
two equivalent sets of defining equations,
\thetag{2.2.9},  for $\scrM$. The
operation $*$ here is taken in 
$\bC^{l}$, i.e. mapping sets in
$\bC^{l}_{(z,w_1,...,w_{l-1})}$ to
$\bC^{l}_ {(\chi,\tau_1,...,\tau_{l-1})}$.)
We claim that the set $S$ is contained in a
proper algebraic variety
$A_l\subset\bC^l$. To see this, note first
that
\thetag{2.3.3} (which, by the induction
hypothesis, holds for
$N_{l-1}$) implies that
$\pi_l(N_{l-1})$ is contained in a proper 
irreducible algebraic variety in
$\bC^{l}$. Let $P_1
(\chi,\tau_1,...,\tau_{l-2})$ be a
(non-trivial) polynomial that vanishes  on
$^*V_{l-1}\subset\bC^{l-1}$, and let
$P_2(\chi,\tau_1,...,\tau_{l-1})$ be a 
(non-trivial) irreducible polynomial that
vanishes on
$^*\pi_l(N_{l-1})$.  Thus, if
$(z,w_1,...,w_{l-1})\in S$
then there exists a
$\chi\in\bC$ such that
$$
\aligned&\tilde
P_1(\chi,z,w_1,...,w_{l-2}):=P_1(\chi,w_1+\bar
q_1(\chi,z), ...,w_{l-2}+\bar
q_{l-2}(\chi,z,w_1,...,w_{l-3}))=0\\&\tilde
P_2(\chi,z,w_1, ...,w_{l-1}):=
P_2(\chi,w_1+\bar q_1(\chi,z),...,w_{l-1}+\bar
q_{l-1}(\chi,z,w_1,
...,w_{l-2}))=0,\endaligned\tag2.3.6
$$ i.e. $\tilde R(z,w_1,...,w_{l-1})=0$ if we
denote by $\tilde R$ the resultant of $\tilde
P_1$ and $\tilde P_2$ as polynomials in
$\chi$. The proof will be complete (with
$A_l=\tilde R^{-1}(0)$) if we can show that
$\tilde R$ is not identically 0, i.e. $\tilde
P_1$ and $\tilde P_2$ have no common factors 
(it is easy to see that neither $\tilde P_1$
nor $\tilde P_2$ is  identically 0). Note
that, for arbitrary
$\tau_1,...,\tau_{l-1}$, we have  (cf.
\thetag{2.2.9})
$$
\tilde
P_2(\chi,z,\tau_1+q_1(z,\chi),...,\tau_{l-1}+q_{l-1}(z,\chi,\tau_1,...,
\tau_{l-2}))=P_2(\chi,\tau_1,...,\tau_{l-1}).\tag2.3.7
$$ It follows from this that
$\tilde P_2$ is irreducible (since
$P_2$ is irreducible). Thus,
$\tilde P_1$ and $\tilde P_2$ cannot have any
common factors because $\tilde P_2$ itself is
the only non-trivial factor of
$\tilde P_2$ and, by the form
\thetag{2.3.3} of $N_{l-1}$,
$\tilde P_2$ is not independent of
$w_{l-1}$. This completes the proof of
Assertion 2.3.3.
\qed\enddemo
 We proceed with the proof of
Proposition 2.3.2. Let us denote by 
$B_l\subset\bC^{l-1}$ the proper algebraic
variety with the property that
$(z,w_1,...,w_{l-2})\in\bC^{l-1}\setminus
B_l$ implies that the polynomial
$\tilde P_1(X,z,w_1,...,w_{l-2})$ defined by
(2.3.6), considered as a polynomial in $X$, has the
maximal number of  distinct roots.  Let
$C_l\subset\bC^l$ denote the union of $A_l$
and
$B_l\times\bC$. For $(z,w_1,...,w_{l-2})$
fixed, let
$\Omega(z,w_1,...,w_{l-2})\subset\bC$ be the
domain obtained by removing from $\bC$ the
roots in $X$ of the polynomial equation 
$$
\tilde P_1(X,z,w_1,...,w_{l-2})=0.\tag2.3.8
$$
 In view of
Assertion 2.3.3 and the inductive hypothesis that 
\thetag{2.3.3} holds for $N_{l-1}$, it
follows from \thetag{2.3.4} that
$$
\aligned&N_l\cap((\bC^l\setminus
C_l)\times\bC^{d-l+1})=\\&
\{(z,w_1,...,w_d)\in\left((\bC^j\setminus
C_l)\times\bC^{d-j+1}\right)
\:\\& \exists
\chi\in\Omega(z,w_1,...,w_{l-2})\subset\bC\,,\,
w_k=g_{lk}(\chi,z,w_1,...,
w_{k-1})\,,\,k=l-1,...,d\},\endaligned
\tag2.3.9
$$ 
where
$$
\aligned g_{lk}(\chi,z,w_1,...,w_{k-1})&=\bar
f_{l-1,k}(\chi,w_1+\bar q_1(\chi,z),...,
w_{l-2}+\bar
q_{l-2}(\chi,z,w_1,...,w_{l-3}))\\&+q_k(z,\chi,w_1+\bar
q_1 (\chi,z), ...,w_{k-1}+\bar
q_{k-1}(\chi,z,w_1,...,w_{k-2})),
\endaligned\tag2.3.10
$$ for $k=l-1,...,d$.  Note that each
$g_{lk}$, for $k=l-1,...,r$, is a
(multi-valued) algebraic  function such that
all branches are holomorphic in a
neighborhood of every point $(\chi,z,w)$
considered in \thetag{2.3.9}, and
$g_{lk}\equiv0$ for $k=r+1,...,d$.

Now, suppose that
$g_{l,l-1}(\chi,z,w_1,...,w_{l-2})$ actually
depends on 
$\chi$, i.e. 
$$
\frac{\partial g_{l,l-1}}{\partial
\chi} (\chi,z,w_1,...,w_{l-2})\not\equiv0.
\tag2.3.11
$$ Then, for each
$(\chi^0,z^0,w_1^0,...,w_{l-2}^0,w_{l-1}^0)$
such that one branch $g$ of
$g_{l,l-1}$ is holomorphic near 
$(\chi^0,z^0,w_1^0,...,w_{l-2}^0)$ with
$$
\frac{\partial g}{\partial \chi}
(\chi^0,z^0,w_1^0,...,w_{l-2}^0)\neq0
\tag2.3.12
$$ and 
$$
w_{l-1}^0=g(\chi^0,z^0,w_1^0,...,w_{l-2}^0),\tag2.3.13
$$ we may apply the (algebraic) implicit
function theorem and deduce that there is a
holomorphic branch $\theta(z,w_1,...,w_{l-1})$
of an algebraic function  near
$(z^0,w_1^0,...,w_{l-1}^0)$ such that
$$
w_{l-1}-g(\theta(z,w_1,...,w_{l-1}),z,w_1,...,w_{l-2})\equiv0.\tag2.3.14
$$ Since $g_{l,l-1}$ is an algebraic
function, which in particular means that any
two choices of branches $g$  at (possibly
different) points
$(\chi^0,z^0,w_1^0,...,w_{l-2}^0)$ can be 
connected via a path in
$(\chi,z,w_1,...,w_{l-2})$ space avoiding the
singularities of
$g_{l,l-1}$ and also avoiding the zeros of
$\partial  g_{l,l-1}/
\partial \chi$, it follows that any solution
$\theta$ of
\thetag{2.3.14} near a point
$(z^0,w_1^0,...,w_{l-1}^0)$ can be
analytically continued to any other solution
near a (possibly  different) point. Thus, all
solutions
$\theta$ are branches of the same algebraic
function, and we denote that algebraic
function by
$\theta_l$. As a consequence, there is an
irreducible polynomial 
$R_l(X,z,w_1,..., w_{l-1})$ such that
$X=\theta_l(z,w_1,...,w_{l-1})$ is its root.
Let
$D_l\subset\bC^{l}$ be the zero locus of the
discriminant of $R_l$ as a polynomial in $X$.
Outside
$(C_l\cup D_l)\times\bC^{d-l+1}
\subset\bC^{d+1}$, we can, by solving for
$\chi=\theta_l(z,w_1,...,w_{l-1})$  in the
equation
$$
w_{l-1}=g_{l,l-1}(\chi,z,w_1,...,w_{l-2}),\tag2.3.15
$$ describe $N_l$ as the (multi-sheeted) graph
$$
w_k=f_{lk}(z,w_1,...,w_{l-1}):=g_{lk}(\theta_l(z,w_1,...,w_{l-1}),w_1,...,
w_{k-1}),\tag2.3.16
$$  for $k=l,...,d$. Clearly, we have
$f_{lk}\equiv0$ for
$k=r+1,...,d$.  By taking $V_l$ to be the
union of $C_l\cup D_l$ and the proper
algebraic variety consisting of points where
any two distinct branches  of $f_{lk}$
coincide (for some $k=l,...,d$), we have
completed the proof of  the inductive step
for $j=l$ under the assumption that $g_{l,l-1}
(\chi,z,w_1,...,w_{l-2})$ actually depends on
$\chi$. 

Now, we complete the proof of the proposition
by showing that Condition 2.3.1 forces
\thetag{2.3.11} to hold as long as 
$l-1\leq r$. Assume, in order to reach a
contradiction, that
$g_{l,l-1}(\chi,z,w_1,..., w_{l-2})$ does not
depend on
$\chi$. It is easy to verify from the form 
\thetag{2.3.1} of $M$ that the sets
$\pi_k(N_j)$, for $j=0,...,k$, are the Segre sets
of $M^k$ at 0. Let us denote these
sets by
$N_j(M^k)$. Now, note that if we pick
$(z^0,w_1^0 ,...,w_{l-1}^0)\in M^l$ then 
$$ (\bar z^0,w_1^0+\bar q_1(\bar
z^0,z^0),...,w_{l-1}^0+\bar q_{l-1}(\bar
z^0,z^0, w_1^0,...,w_{l-2}^0))=(\bar z^0,\bar
w_1^0,...,\bar w_{l-1}^0).\tag2.3.17
$$ Thus, if we pick the point
$(z^0,w_1^0,...,w_{l-1}^0)\in M^l$ such that 
it is not on the algebraic variety $C_l$
(which is possible since the generic real
submanifold $M^l$ cannot be contained in a
proper algebraic variety; $C_l\cap M^l$ is a
proper real algebraic subset of $M^l$) then,
by construction of $C_l$, the point
$$ (\bar z^0,w_1^0+\bar q_1(\bar
z^0,z^0),...,w_{l-2}^0+\bar q_{l-2}(\bar
z^0,z^0, w_1^0,...,w_{l-3}^0))=(\bar z^0,\bar
w_1^0,...,\bar w_{l-2}^0)\tag2.3.18
$$ is not in $^*\pi_l(V_{l-1})$. By the
induction hypothesis,
$\pi_l(N_{l-1})=N_{l-1}(M^l)$ consists of a
$b_{l-1,l-1}$-sheeted graph (each sheet,
disjoint from the other, corresponds to a
branch of 
$f_{l-1,l-1}$) above a neighborhood of the
point $(z^0,w_1^0,...,$
$w_{l-2}^0)$. Since $g_{l,l-1}$ is assumed
independent of $\chi$, we can, in view of
\thetag{2.3.18}, take $\chi=\bar z$ in the
defining equation
$$
w_{l-1}=g_{l,l-1}(\chi,z,w_1,...,w_{l-2})\tag2.3.19
$$ for $N_l(M^l)$, near the point
$(z^0,w_1^0,...,w_{l-2}^0)$.  {}From the
definition \thetag{2.3.10} of
$g_{l,l-1}$ and \thetag{2.3.18} it follows
that $N_l(M^l)$ also consists of a
$b$-sheeted graph, with
$b\leq b_{l-1,l-1}$, (each sheet
corresponds to a choice of branch of 
$\bar f_{l-1,l-1}$ at $(\bar z^0,\bar w_1^0,$
$...,\bar w_{l-2}^0)$)  above a neighborhood
of the point
$(z^0,w_1^0,...,w_{l-2}^0)$. Since
$N_{l-1}(M^l)\subset N_{l}(M^l)$, we must
have $b=b_{l-1,l-1}$ and, moreover, for each
branch
$f^k_{l-1,l-1}$ there is  possibly another
branch
$f^{k^\prime}_{l-1,l-1}$ such that for every
$(z,w_1,...,w_{l-2})$ the following holds
$$
\aligned f^{k}_{l-1,l-1}(z,w_1,...,
&w_{l-2})=\\&
\bar f^{k^\prime}_{l-1,l-1}(\bar z,w_1+\bar
q_1(\bar z,z),..., w_{l-2}+\bar q_{l-2}(\bar
z,z,w_1,...,w_{l-3}))\\&+q_{l-1}(z,\bar
z,w_1+\bar  q_1(\bar z,z),...,w_{l-2}+\bar
q_{l-2}(\bar z,z,w_1,...,w_{l-3})).
\endaligned\tag2.3.20
$$ Since all the sheets of the graphs are
disjoint, the mapping $k
\to k^\prime$ is a permutation. We average
over $k$ and $k^\prime$, restrict to points
$(z,w_1,...,w_{l-2})\in M^{l-1}$, and obtain,
by \thetag{2.3.18} and
\thetag {2.3.20},
$$
\aligned
\frac{1}{b_{l-1,l-1}}\sum_{k=1}^{b_{l-1,l-1}}
f^{k}_{l-1,l-1}(z,w_1,...,w_{l-2}) &=
\frac{1}{b_{l-1,l-1}}\sum_{k^\prime=1}^{b_{l-1,l-1}}
\bar f^{k^\prime}_{l-1,l-1}(\bar z,\bar
w_1,...,
\bar w_{l-2})\\&+q_{l-1}(z,\bar z,\bar
w_1,...,\bar w_{l-2}).
\endaligned\tag2.3.21
$$ Let us denote by $f$ the holomorphic
function near
$(z^0,w_1^0,...,w_{l-2}^0)$ defined by 
$$
f(z,w_1,...,w_{l-2})=\frac{1}{b_{l-1,l-1}}\sum_{k=1}^{b_{l-1,l-1}}
f^{k}_{l-1,l-1}(z,w_1,...,w_{l-2}),\tag2.3.22
$$ and by $K\subset\bC^l$ the CR manifold of
CR dimension 1 defined near
$(z^0,w_1^0,...,w_{l-2}^0,$
$f(z^0,w_1^0,...,w_{l-2}^0))$ by
$$
K:=\{(z,w_1,...,w_{l-1})\:(z,w_1,...,w_{l-2})\in
M^{l-1}\,,\,w_{l-1}=
f(z,w_1,...,w_{l-2})\}.\tag2.3.23
$$ The equation \thetag{2.3.21} immediately
implies that $K\subset M^l$. By Condition
2.3.1, $M^l$ is of finite type near 0. Note
that, by the form
\thetag{2.3.1} of $M$, the condition that
$M^l$ is of  finite type at a point  is only a
condition on $(z,w_1,...,w_{l-2})$ (i.e. not
on $w_{l-1}$). Thus, by picking the point
$(z^0,w_1^0,...,w_{l-2}^0)\in M^{l-1}$
sufficiently close to 0 (which is possible
since, as we mentioned above, $C_l\cap M^l$
is a proper real algebraic subset of
$M^l$), we reach the desired contradiction.
This completes the proof of Proposition
2.3.2.\qed\enddemo 

\subhead 2.4. Homogeneous submanifolds of
arbitrary CR dimension\endsubhead We prove
here the analog of Proposition 2.3.2 for
arbitrary CR dimension $n$. 
\proclaim{Proposition 2.4.1} Let
$M$ be of the form \thetag{2.3.1} and assume
that $M$ satisfies Condition {\rm 2.3.1}. Let
$N_0,N_1,...,N_{j_0}$ be the Segre sets
 of $M$ at $0$. Then, for
each $j=1,...,j_0$, there is a partition of
the set
$\{1,2,...,r\}$ into
$I_j=\{i_1,i_2,...,i_{a_j}\}$ and
$K_j=\{k_1,k_2,...,k_{b_j}\}$ such that
$$
\emptyset=I_1\ssneq I_2\ssneq
I_3\ssneq...\ssneq I_{j_0}=
\{1,2,...,r\},\tag2.4.1
$$ and there is a proper algebraic variety
$V_j\subset\bC^{n+a_j}$ such that 
$N_j$ satisfies 
$$
\aligned&N_j\cap\left((\bC^{n+a_j}\setminus V_j)\times\bC^{b_j}\times
\bC^{d-r}\right)=\\
&\left\{(z,w_1,...,w_d)\: \left\{
\aligned w_{k_\mu} &=f_{jk_{\mu}}(z,w_{i_1},...,w_{i_{a_j}})\,,\,\mu=1,...,b_j
\\w_k &=0\hskip 35truemm,\,k=r+1,...,d\endaligned\right.\right\}.
\endaligned\tag2.4.2
$$
 Here  $(z,w_{i_1},..., w_{i_{a_j}}) \in
\bC^{n+a_j}$ and
$(w_{k_1},...,w_{k_{b_j}}) \in \bC^{b_j}$.
Each
$f_{jk_{\mu}}$, for
$k=1,...,b_j$,  is a  {\rm(}multi-valued{\rm)}
algebraic  function with
$b_{jk_{\mu}}$ holomorphic, disjoint
bran\-ches outside $V_j$ and, moreover, each
$f_{jk_{\mu}}$ is independent of $w_{i_\nu}$
for all
$i_\nu>k_\mu$.\endproclaim
\demo{Proof}  We emphasize here those aspects of the
proof which are different from that of
Proposition 2.3.2. We proceed by induction  on
$j$. The statement of the
proposition holds for
$j=1$ with
$V_1=\emptyset$ and each
$f_{1k_{\mu}}\equiv0$. Assume the statement
holds for $j=1,...,l-1$. Let us for
simplicity denote the numbers $a_{l-1}$ and
$b_{l-1}$ by
$a$ and
$b$, respectively. The representation
\thetag{2.3.4}, with
$z$ in $\bC^n$ rather than $\bC$, still
holds. Let
$\Omega(z,w)\subset\bC^n$ be the complement
of the algebraic subset of
$\chi$ such that, for fixed
$(z,w)\in\bC^{n+d}$, 
$$ (\chi,w_{i_1}+
\bar q_{i_1}(\chi,z,w_1,...,w_{i_1-1}),...,
w_{i_{a}}+\bar
q_{i_{a}}(\chi,z,w_1,...,w_{i_{a}-1}))\in
{}^*V_{l-1}.\tag2.4.3
$$ We describe a part $\tilde N_l$ of
$N_{l}$ as follows
$$
\aligned\tilde
N_l=\{&(z,w)\in\bC^{n+d}\:\exists
\chi\in
\Omega(z,w)\subset\bC^n\,,\,\\&
w_{k_\mu}=g_{lk_{\mu}}(\chi,z,w_1,...,w_{k_{\mu}-1})\,,\,\mu=1,...,b\,;\,
w_k=0\,,\,k=r+1,...,d\},
\endaligned\tag2.4.4
$$ 
where 
$$
\aligned
g_{lk_{\mu}}&(\chi,z,w_1,...,w_{k_\mu-1})=\\&\bar
f_{l-1,k_\mu}(\chi,w_{i_1}+
\bar q_{i_1}(\chi,z,w_1,...,w_{i_1-1}),...,
w_{i_{a}}+\bar
q_{i_{a}}(\chi,z,w_1,...,w_{i_{a}-1}))\\
&+q_{k_\mu}(z,\chi,w_1+\bar q_1 (\chi,z),
...,w_{k_\mu-1}+\bar
q_{k_\mu-1}(\chi,z,w_1,...,w_{k_\mu-2}))
\endaligned\tag2.4.5
$$ The fact that $w_k=0$ for
$k=r+1,...,d$ follows from
\thetag{2.4.2} with 
$j=l-1$ and the form \thetag{2.3.1} of
$M$. Note also that, by the induction
hypothesis,
$f_{l-1,k_\mu}$ is independent of
$w_{i_\nu}$ for $i_\nu>k_\mu$. Let
$w^\prime=(w_{i_1},...,w_{i_a})$ and
$w^{\prime\prime}= (w_{k_1},...,w_{k_b})$.
Note that , for generic 
$(z,w^{\prime\prime})\in\bC^{ n+b}$, the
mapping from $\bC^{n+a}$ into itself given by
$$
(\chi,w^\prime)\mapsto(\chi,w_{i_1}+
\bar q_{i_1}(\chi,z,w_1,...,w_{i_1-1}),...,
w_{i_{a}}+\bar
q_{i_{a}}(\chi,z,w_1,...,w_{i_{a}-1}))\tag2.4.6
$$ has generic rank $n+a$ (indeed, it has
rank $n+a$ near the origin for $z=0$). Thus,
the set of
$w^\prime\in\bC^{a}$ for which
\thetag{2.4.3} holds (with small
$z$ and $w^{\prime\prime}$  arbitrary) for
all $\chi\in\bC^n$ is a proper algebraic
variety.  Restricting $(\chi,w^\prime)$ to
the complement of the set where
\thetag{2.4.3} holds, we consider the mapping
\thetag{2.4.6} with 
$$
w_{k_1}=g_{lk_1}(\chi,z,w_1,...,w_{k_1-1})\tag2.4.7
$$ instead of $w_{k_1}$ fixed. Again, one
verifies that this mapping has generic rank
$n+a$ for generic $(z,w_{k_2},...,w_{k_{b}})$
(e.g. with $z$ small), and thus the set of
$w^\prime$ for which \thetag {2.4.3} holds
(with $w_{k_1}$ given by \thetag{2.4.7}) for
all
$\chi$ is a proper algebraic variety. By
proceeding inductively, substituting
$g_{lk_\mu}$ for $w_{k_\mu}$ in the mapping
\thetag{2.4.6}, we find that we can take for
$\tilde N_l$ (for brevity, we write
$w^{\prime\prime\prime}= (w_{r+1},...,w_d)$)
$$
\aligned\tilde
N_l=\{&(z,w^\prime,w^{\prime\prime},w^{\prime\prime\prime})
\in(\bC^{n+a}\setminus
C_l)\times\bC^{b}\times\bC^{d-r}\:\exists
\chi\in\\&
\tilde\Omega(z,w^\prime)\subset\bC^n\,,\,
w_{k_\mu}=\tilde
g_{lk_{\mu}}(\chi,z,w^\prime)\,,\,\mu=1,...,b\,;\,w^
{\prime\prime\prime}=0\},
\endaligned\tag2.4.8
$$ where $C_l\subset\bC^{n+a}$ is a proper
algebraic  variety, $\tilde
\Omega(z,w^\prime)\subset\bC^n$ is the
complement of a proper algebraic variety in
$\bC^n$, and where $\tilde g_{lk_1}= g_{lk_1}$
and subsequent $\tilde g_{lk_\mu}$ are
obtained from $g_{lk_\mu}$  by substituting
$$ w_{k_\gamma}=\tilde
g_{lk_\gamma}(\chi,z,w^\prime)\quad,
\quad\gamma=1,...,\mu-1.\tag2.4.9
$$ Thus, each $\tilde g_{lk_\mu}$ is a
function only of those
$w_{i_1},..., w_{i_\nu}$ for which
$i_\nu<k_\mu$. 

As in the proof of Proposition
2.3.2, we  assume first
that the map 
$$
\bC^{2n+a}\ni(\chi,z,w^\prime)\mapsto (\tilde
g_{lk_1} (\chi,z,w^\prime),...,\tilde
g_{lk_b}(\chi,z,w^\prime))=:G(\chi,z,w^\prime)
\in\bC^b\tag2.4.10
$$ actually depends on $\chi$, i.e. 
$$ G_\chi(\chi,z,w^\prime):=\frac{\partial
G}{\partial \chi}(\chi,z,w^\prime)
\not\equiv0.\tag2.4.11
$$ Denote by $m\geq1$ the maximal rank
of $G_\chi$, and by
$G^\prime=(G^{t_1},...,G^{t_m})$ the $m$
first components of $G$ such  that
$G^\prime_\chi$ has generic rank $m$ (thus,
the set
$\{t_1,...,t_m\}$ is a subset of
$K_{l-1}$). Note that this does not
necessarily need to be the first
$m$ components of $G$, but any component
$G^t$, with
$t_\alpha<t<t_{
\alpha+1}$ for some
$\alpha\in\{1,...,m-1\}$, has then the
property that 
$$
G^t_\chi(\chi,z,w^\prime)\equiv\sum_{j=1}^\alpha
c_j(\chi,z,w^\prime)
G^{t_j}_\chi(\chi,z,w^\prime),\tag2.4.12
$$ for some functions
$c_1,...,c_\alpha$. We may assume, by an
algebraic change of  coordinates in the
$\chi$ space if necessary, that
$G^\prime_{\chi^\prime}$, where
$\chi^\prime=(\chi_1,...,\chi_m)$, has
generic rank $m$ and that $G$ is independent
of the last coordinates $\chi^{\prime\prime}:=
(\chi_{m+1},...,\chi_n)$. Now, solve for
$\chi^\prime=\theta_l(z,w^\prime,w_{t_1},...,w_{t_m})$
in the equations
$$ w_{t_j}=G^{t_j}(\chi,z,w^\prime)\quad,\quad
j=1,...,m.\tag2.4.13
$$ The solution
$\theta_l$ is a (multi-valued) algebraic
function. By substituting
$$
\chi^\prime=\theta_l(z,w^\prime,w_{t_1},...,w_{t_m})
$$ in the remaining equations for
$\tilde N_l$ (and remembering that, by the
choice of $m$ and 
$\chi^\prime$, these equations are
independent of
$\chi^{\prime\prime}$) we find, denoting  by
$K_l:=\{u_1,...,u_{b-m}\}$ the complement of
the set
$\{t_1,...,t_m\}$ in $K_{l-1}$,
$$
\aligned w_{u_j} &=\tilde
g_{lu_j}(\theta_l(z,w^\prime,w_{t_1},...,w_{t_m}),\chi^
{\prime\prime},z,w^\prime)\\&=:f_{lu_j}(z,w^\prime,w_{t_1},...,w_{t_m})\quad,
\quad j=1,...,b-m.\endaligned\tag 2.4.14
$$ Since $\tilde N_l$ is a dense open subset
of
$N_l$, the equations
\thetag{2.4.14} imply that $N_l$ is indeed
of the form
\thetag{2.4.2}, with $K_l\ssneq K_{l-1}$ as
defined above and
$I_l=\{1,...,r\}\setminus K_l$, and where we
let
$V_l\subset\bC^{n+a+m}$ be a suitable proper
algebraic variety containing the
singularities of the algebraic functions
$f_{lu_j}$ ($j=1,...,b-m$). To finish the
proof (under the assumption that the mapping
$G$ actually depends on $\chi$), we need to
show that each $f_{lu_\nu}$ is independent of
$w_{t}$ for
$t>u_\nu$. Recall that
$$
f_{lu_\nu}(z,w^\prime,w_{t_1},...,w_{t_m})=G^{u_\nu}(\theta_l(z,w^\prime,
w_{t_1},...,w_{t_m}),\chi^{\prime\prime},z,w^\prime).\tag2.4.15
$$ Let $1\leq\alpha<m-1$ be the number such
that
$t_\alpha<u_\nu<t_{\alpha+1}$ (unless there
is such a number there is nothing to prove),
and differentiate
\thetag{2.4.15} with respect to
$w_{t}$, where $t\geq u_\nu$. Using
\thetag{2.4.12}, we obtain (using vector
notation; recall that
$G^{u_\nu} (\xi,z,w^\prime)$ is independent
of $w_t$)
$$
 f_{lu_\nu,w_{t}}
=G^{u_\nu}_{\chi^\prime}\theta_{l,w_{t}}
=\sum_{j=1}^\alpha
c_jG^{t_j}_{\chi^\prime}\theta_{l,w_{t}}.
\tag2.4.16
$$ Now, by the definition of
$\theta_l$,
$
G^{t_j}(\theta_l(z,w^\prime,w_{t_1},...,w_{t_m}),\chi^{\prime\prime},z,
w^\prime)\equiv w_{t_j},
$ and so 
$
G^{t_j}_{\chi^\prime}\theta_{l,w_{t}}=0,\quad\text{\rm
if }t>t_j.
$ Thus, since $t>t_j$ for
$j=1,...,\alpha$, it follows from
\thetag{2.4.16} that
$ f_{lu_\nu,w_{t}}=0.
$ This proves the induction hypothesis for
$j=l$ under the assumption that the mapping
$G$ actually depends on $\chi$. 

As in the proof of Proposition 2.3.2, we
are left to show that Condition
2.3.1 implies that $G$ actually depends on 
$\chi$ as long as
$I_{l-1}\ssneq\{1,2,...,r\}$. Assume,  in
order to reach a contradiction that $G$ does
not depend on $\chi$. In particular then, the
function $\tilde g_{lk_1}= g_{lk_1}$ does not
depend on $\chi$. Since, by the induction
hypothesis,
$g_{lk_1}(\chi,z,w^\prime)$ does not depend
on $w_j$ for $j\geq k_1$, we can consider the
projection $\pi_{k_1}$ and proceed exactly as
in the conclusion of the proof of Proposition
2.3.2. We leave the straightforward
verification to the reader. The proof
of Proposition 2.4.1 is now
complete.\qed\enddemo 

\subhead 2.5 Proof of Theorem
2.2.1\endsubhead By the remarks preceding
the theorem, we may assume $M$ is generic
throughout this proof.  We start by
proving (a). Since the Segre sets of
$M$ at $p_0$ are invariantly defined, we may
choose any holomorphic coordinates
near
$p_0$. Let $m_1\le...\le m_r$ be the 
H\"ormander numbers of $M$ at $\po$.  By
[BR1, Theorem 2], there are holomorphic
coordinates
$(z,w)\in\bC^n\times\bC^d$ such that the
equations of $M$ near
$p_0$ are given by
$$
\left\{\aligned w_j &=\bar w_j+q_j(z,\bar
z,\bar w_1,...,\bar w_{j-1})+R_j (z,\bar
z,\bar w)
\quad,\quad j=1,...,r\\w_k  &=\bar
w_k+\sum_{l=r+1}^d f_{kl}(z,\bar z,\bar w)
\bar w_l\quad,\quad
k=r+1,...,d,\endaligned\right.\tag2.5.1
$$ where, for $j=1,...,r$,
$q_j(z,\bar z,\bar w_1,...,\bar w_{j-1})$ is
weighted  homogeneous of degree $m_j$,
$R_j(z,\bar z,\bar w)$ is a real analytic
function whose Taylor expansion at the
origin consists of terms  of weights at least
$m_j+1$, and the
$f_{kl}$ are real analytic functions that
vanish at the origin.  Here, $z$ is assigned
the weight 1, $w_j$ the weight $m_j$ for
$j=1,...,r$ and weight $m_r+1$ for
$j=r+1,...,d$. Moreover, 
 the homogeneous
manifold $M^0\subset \cnn$ defined by
$$
\left\{\aligned w_j &=\bar w_j+q_j(z,\bar
z,\bar w_1,...,\bar w_{j-1})
\quad,\quad j=1,...,r\\w_k  &=\bar
w_k\quad,\quad k=r+1,...,d\endaligned
\right.\tag2.5.2
$$ satisfies Condition 2.3.1. For $\epsilon
> 0$,  we  introduce the scaled coordinates
$(\tilde z,\tilde w)\in\bC^{n+d}$ defined by
$$
\left\{\aligned z &=z(\tilde
z;\epsilon)=\epsilon\tilde z\\w_j &=
w_j(\tilde w;\epsilon)=\epsilon^{l_j}\tilde
w_j\quad,\quad
j=1,...,d,\endaligned\right.\tag2.5.3
$$ where $l_j=m_j$ for $j=1,...,r$ and
$l_k=m_r+1$ for $k=r+1,...,d$. We write
$\tilde f_{kl}$ for the function
$$
\tilde f_{kl}(\tilde z,\bar{\tilde
z},\bar{\tilde w};\epsilon)=
\frac{1}{\epsilon}f_{kl}(z(\tilde
z;\epsilon),\bar z(\tilde z;\epsilon),
\bar w(\tilde w;\epsilon)),\tag2.5.4
$$ and similarly,
$$
\tilde R_{j}(\tilde z,\bar {\tilde z},\bar
{\tilde w};\epsilon)=
\frac{1}{\epsilon^{m_j+1}}R_{j}(z(\tilde
z;\epsilon),\bar z(\tilde z;\epsilon),
\bar w(\tilde w;\epsilon)).\tag2.5.5
$$ Note that both $\tilde f_{kl}(\tilde
z,\bar{\tilde z},\bar{\tilde w};\epsilon)$ and
$\tilde R_j(\tilde z,\bar{\tilde
z},\bar{\tilde w};\epsilon)$ are real
analytic functions of
$(\tilde z,\tilde w;\epsilon)$ in a
neighborhood of $(0,0;0)$. In the scaled 
coordinates,
$M$ is represented by the equations
$$
\left\{\aligned \tilde w_j &=\bar{\tilde
w}_j+q_j(\tilde z,\bar{\tilde z},
\bar  {\tilde w}_1,...,\bar{\tilde
w}_{j-1})+\epsilon\tilde R_j (\tilde
z,\bar{\tilde z},\bar{\tilde w};\epsilon)
\quad,\quad j=1,...,r\\\tilde w_k 
&=\bar{\tilde w}_k+\epsilon
\sum_{l=r+1}^d \tilde f_{kl}(\tilde
z,\bar{\tilde z},\bar{\tilde w};\epsilon)
\bar {\tilde w}_l\quad,\quad
k=r+1,...,d,\endaligned\right.\tag2.5.6
$$ Now, let $\tilde
v^j(\tilde\Lambda;\epsilon)$ be the mapping
$\bC^{jn}\mapsto
\bC^d$, described in \S2.2, such that the
Segre set $N_j$ of $M$ at
$p_0$ is parametrized by 
$$
\bC^{jn}\ni\tilde\Lambda\mapsto (\tilde
z,\tilde v^j(\tilde\Lambda;
\epsilon))\in\bC^N\tag2.5.7
$$  in the scaled coordinates
$(\tilde z,\tilde w)$ (cf. 
\thetag{2.2.10}--\thetag{2.2.12} and
\thetag{2.2.13}--\thetag{2.2.15}  to see how
the map \thetag{2.5.7} is obtained from the
defining equations
\thetag{2.5.6}). Note that 
$\tilde v^j$ depends real analytically on the
 small parameter
$\epsilon$. The
generic dimension $d_j$ of the Segre set 
$N_j$ is the generic rank of the mapping
\thetag{2.5.7} with
$\epsilon\neq0$, and is in fact independent
of $\epsilon $. By the real  analytic
dependence on
$\epsilon$ there is a neighborhood $I$ of
$\epsilon=0$ such that the generic rank of
\thetag{2.5.7}, for all
$ \epsilon \in I\backslash \{0\}$, is at least
the generic rank of
\thetag{2.5.7} with $\epsilon=0$. For
$\epsilon=0$  the mappings
\thetag{2.5.7} parametrize the Segre sets
$N^0_j$ of the homogeneous manifold
$M^0$ defined by 
\thetag{2.5.2}. By Proposition
2.4.1, applied to the Segre sets $N^0_j$ of
$M^0$ at 0, we deduce that the generic dimension
of the maximal Segre set of $M^0$ at 0 is
$n+r$. Thus, $d_{j_0}\ge n+r$, where
$d_{j_0}$ is the generic dimension of the
maximal Segre set of $M$ at $\po$ . On the
other hand, if we  go back to the unscaled
cordinates 
$(z,w)$, we note from the construction of
the Segre sets
 that each $N_j$ is contained
in the complex manifold
$X=\{(z,w)\:w_{r+1}=...=w_d=0\}$. Thus
$d_{j_0} \le n+r$, so that we obtain the
desired equality $d_{j_0} = n+r$.  This
proves part (a) of the theorem. 

It follows from \thetag{2.5.1} that the CR
vector fields of $M$
are all tangent to
$M\cap X=\{(z,w)\in
M\:w_j=0\,,\,j=r+1,...,d\}$. Thus, the local
CR orbit of $p_0$ is contained in $M\cap X$.
Also, since there are $r$  H\"ormander
numbers, the CR orbit of $p_0$ has
dimension $2n+r$. Since the dimension of
$M\cap X$ is $2n+r$ as well, it follows that
the local CR orbit of $p_0$ is $M\cap X$.
This proves part (b) of the  theorem. 

Finally, to prove part (c) of the theorem we
note that if $M$ is real algebraic then each
Segre set $N_j$ is contained in a unique
irreducible complex algebraic variety of
dimension $d_j$. Since
$N_{j_0}$ contains a relatively open subset
of $X$, this relatively open subset of $X$
coincides with a relatively open subset of the
unique algebraic variety containing
$N_{j_0}$. Hence, $X$ is complex algebraic.
This completes the proof of  Theorem
2.2.1.\qed

\heading 3. Algebraic properties of holomorphic
mappings between real algebraic sets\endheading
\subhead 3.1. A generalization of
Theorems 1 and 4\endsubhead 
We denote by
$\scrO_N (p_0)$ the ring of germs of holomorphic
functions in $\bC^N$ at $p_0$, and by 
$\scrA_N(p_0)$ the subring of
$\scrO_N(p_0)$ consisting of those germs that are
also algebraic, i.e. those germs for which there is a
nontrivial polynomial
$P(Z,x)\in\bC[Z,x]$ (with
$Z\in\bC^N$ and $x\in\bC$) such that any
representative $f(Z)$ of the germ satisfies
$$ P(Z,f(Z))\equiv0.\tag3.1.1
$$ In particular, any function in
$\scrA_N(p_0)$ extends as a possibly multi-valued
holomorphic function in $\bC^N\setminus V$, where $V$
is a proper algebraic variety in $\bC^N$.  We refer
the reader to e.g. [BR3, \S1] for some elementary
properties of algebraic holomorphic functions that
will be used in this paper.  If $U\subset\bC^N$ is a
domain we denote by
$\scrO_N(U)$ the space of holomorphic functions
in $U$.

If $X\subset\cnn$ is an algebraic variety with
$\dim\ X=K$,
$\po
\in X_{\text{\rm reg}}$, and
$f$ is a holomorphic function on
$X$ defined near $p_0$ then we say that $f$ is
algebraic if, given algebraic coordinates 
$$
\bC^K\ni\zeta\mapsto Z(\zeta)\in\bC^N\tag3.1.2
$$ on $X$ near $p_0$ with
$Z(0)=p_0$ (i.e. each component of
\thetag{3.1.2} is  in
$\scrA_K(0)$), the function
$h=f\circ Z$ is in $\scrA_K(0)$. The transitivity
property of algebraic functions (e.g. [BM] or
[BR3, Lemma 1.8 (iii)]) implies that this
definition is independent of the choice of
algebraic coordinates. If
$f$ is algebraic on $X$ near $p_0$ and $X$ is
irreducible, then $f$ extends as a possibly
multi-valued holomorphic function on
$X_{\text{reg}}\setminus V$, where
$V$ is a proper algebraic subvariety of
$X_{\text{reg}}$.  (Note that  
$X_{\text{reg}}$ is a connected manifold.) We denote
by $\scrO_X(p_0)$ the ring of germs of holomorphic
functions on $X$ at $p_0$, and by
$\scrA_X(p_0)$ the subring of germs that are
algebraic.

Also, given two points
$p_0\in\bC^N$ and
$p_0^\prime\in\bC^{N^\prime}$, we denote by
$\hol(p_0,p_0^\prime)$ the space of germs of
holomorphic mappings at $p_0$ from $\bC^N$ into
$\bC^{N^\prime}$ taking $p_0$ to
$p_0^\prime$. We denote by 
$\alg(p_0,p_0^\prime)$ the subspace of
$\hol(p_0,p_0^\prime)$ consisting of those germs for
which each component of the mapping is algebraic.
Similarly, given an algebraic variety $X$ in
$\bC^{N}$ with $p_0\in X_{\text{reg}}$, we denote by
$\hol_X(p_0,p_0^\prime)$ the space of germs of 
holomorphic mappings at $p_0$ from $X$ into
$\bC^{N^\prime}$ taking $p_0$ to
$p_0^\prime$, and by
$\alg_X(p_0,p_0^\prime)$ the subspace of germs with
algebraic components. 

Before we present the main results, we state the
following lemma, whose proof is straightforward
and left to the reader.
\proclaim{Lemma 3.1.1} Let
$M$ be a generic real analytic submanifold in 
$\bC^K$ and let $p_0\in M$. Suppose that there is
$h=(h_1,...,h_q)\in (\scrO_K(p_0))^q$ such that the
following holds.
\roster
\item"(i)" $h(p_0)=0$ and $\partial
h_1\wedge...\wedge\partial h_q\neq0$ in a
neighborhood of
$p_0$.
\item"(ii)" $h|_M$ is valued in
$\bR^q$. 
\endroster Then $M\cap S_0$, where
$S_0=\{Z\:h(Z)=0\}$, is a generic submanifold of 
$S_0$ near $p_0$.\endproclaim

We are now in a position to formulate one of the main
results in this paper.

\proclaim{Theorem 3.1.2} Let $M$ be a real algebraic,
holomorphically non-degenerate, CR submanifold in
$\bC^N$, let $\scrV \subset \cnn$ be the smallest
complex algebraic variety containing $M$, and let
$p_0\in\overline{M}$ be a regular point of $\scrV$. 
Assume that there is
$h=(h_1,...,h_q)\in (\scrA_\scrV(p_0))^q$ satisfying
the following.
\roster
\item"(i)" $h(p_0)=0$ and $\partial_\scrV
h_1\wedge...\wedge\partial_\scrV h_q\neq0$ in a
neighborhood of
$p_0$.
\item"(ii)" $h|_M$ is valued in
$\bR^q$. 
\endroster Let $U$ be a sufficiently small
neighborhood of $\po$ in $\cnn$ and denote by $S_c$,
for
$c\in\bC^q$ with 
$|c|$ small,   the algebraic manifolds
$$ S_c=\{Z\in \scrV\cap U\:h(Z)=c\}.\tag3.1.3
$$ Assume that the generic submanifold $M\cap S_{h(p)}$
is minimal at
$p$ for some
$p\in M\cap U$. 

Then if
$A^\prime$ is a real  algebraic set in 
$\bC^{N^\prime}$ with $\dim_\bR
A^\prime=\dim_\bR M$,
$p_0^\prime\in A^\prime$, and  
$H\in\hol_\scrV(p_0,p_0^\prime)$ satisfies
$H(M)\subset A^\prime$, with generic rank equal to
$\dim_\bC\scrV$, there
exists $\delta>0$ such that $H|_{S_c}$ is algebraic
for every $|c|<\delta$.\endproclaim

Note that $M$ is not required to be closed in Theorem
3.1.2. Since
$M$ is real algebraic, it is contained in a real
algebraic set
$A$ of the same dimension in
$\bC^N$ such that $A$, in turn, is contained in the
complex algebraic variety $\scrV$. Thus, the point
$p_0\in \overline {M}$ is a point on the real
algebraic set $A$, and the only thing required of
$p_0$ is that it is a smooth point of
$\scrV$; if e.g.
$M$ is generic then, of course,
$\scrV$ is the whole space $\bC^N$ and, hence,
nothing at all is required of $p_0\in\overline{M}$.
The point $p_0$ could be a singular point of $A$, a
regular point where the CR dimension increases, or 
a point across which $M$ extends as a CR manifold. 

Specializing Theorem 3.1.2 to the case $q=0$ we
obtain the following result.

\proclaim{Corollary 3.1.3} Let $M$ be a real
algebraic, holomorphically nondegenerate, CR
submanifold in $\bC^N$, let $\scrV$ be the smallest
complex algebraic variety that contains $M$, and let
$p_0\in\overline{M}$ be a regular point of $\scrV$. 
Assume that there exists $p\in
M$, such that $M$ is minimal at $p$. 

Suppose
$A^\prime$ is a real  algebraic set in 
$\bC^{N^\prime}$ such that $\dim_\bR
A^\prime=\dim_\bR M$,
$p_0^\prime\in A^\prime$, and  
$H\in\hol_\scrV(p_0,p_0^\prime)$ satisfying
$H(M)\subset A^\prime$, and with generic rank equal to
$\dim_\bC\scrV$. Then
$H\in\alg_\scrV(p_0,p_0^\prime)$. \endproclaim

Specializing again in Corollary 3.1.3 to the case
where $M$ is generic, we obtain the following.

\proclaim {Corollary 3.1.4} Let $M \subset \cnn$ be a
real algebraic, holomorphically nondegenerate,
generic submanifold.
Assume there exists $p\in
M$, such that $M$ is minimal at $p$. 

Suppose
$A^\prime$ is a real  algebraic set in 
$\bC^{N^\prime}$ such that $\dim_\bR
A^\prime=\dim_\bR M$,
$p_0^\prime\in A^\prime$, and  
$H\in\hol(p_0,p_0^\prime)$ satisfying
$H(M)\subset A^\prime$, with generic rank equal to
$N$. Then
$H$ is algebraic.\endproclaim
 
\noindent{\bf Example 3.1.5.} Consider  the generic
holomorphically nondegenerate submanifold
$M\subset\bC^4$ given by
$$
\left\{\aligned \im w_1 &=|z|^2+\re w_2|z|^2\\
\im w_2 &=\re w_3|z|^4\\\im w_3 &=0.
\endaligned\right.\tag3.1.4
$$ The function $h_1(z,w)=w_3$ is real on $M$, and
$M\cap\{(z,w)\:w_3=c\}$ is clearly minimal near
$(z,w_1,w_2)=(0,0,0)$ for all real $c\neq0$. Thus,
Theorem 3.1.2 implies that any holomorphic mapping
$H\:\bC^4\mapsto
\bC^{N^\prime}$ near 0, generically of rank 4, such
that $H(M)$ is contained in a 5 dimensional real
algebraic subset of
$\bC^{N^\prime}$ is algebraic on the leaves
$\{w_3=c\}$, for all sufficiently small
$c\in\bC$.  This result is optimal, because it is
easy to verify that the mapping
$H\:\bC^4\mapsto\bC^4$, defined by
$$ H(z,w_1,w_2,w_3)=(ze^{iw_3},w_1,w_2,w_3),
$$ is a biholomorphism near the origin,  and maps $M$
into itself.  Moreover, $H$ is algebraic on each
$\{w_3=c\}$ but not in the whole space $\bC^4$.

It is interesting to note that the only H\"ormander
numbers at 0 is
$2$, and  that the maximal Segre set
of $M$ at 0 is
$ N_2=\{(z,w)\:w_2=w_3=0\}
$ Thus, the dimension of the maximal Segre set at 0
is smaller than the  dimension of the leaves on which
$H$ is algebraic. For generic points
$p\in M$ though, the maximal Segre set of $M$ at $p$
coincides with one of these leaves. 

\medskip
\noindent{\bf Example 3.1.6.} Consider the real
algebraic subset
$A\subset
\bC^4$ defined by
$$
\left\{\aligned (\im w_1)^2 &=\re
w_2(|z_1|^2+|z_2|^2)\\\im w_2 &=0.
\endaligned\right.\tag3.1.5
$$ It is singular on
$\{(z_1,z_2,w_1,w_2) =
(0,0,s_1,s_2)\:s_1,s_2\in\bR\}$, but outside that set
it is a  generic holomorphically non-degenerate
manifold
$M$. The function
$h_1(z,w)=w_2$ is real on $M$, and
$M\cap\{w_2=s_2^0\}$ is minimal everywhere for all
real
$s_2^0\neq0$.  Theorem 3.1.2 implies that any
holomorphic mapping $H\:\bC^4\mapsto
\bC^{N^\prime}$ near 0, generically of rank 4, such
that $H(M)$ is contained in a 6 dimensional real
algebraic subset of
$\bC^{N^\prime}$ is algebraic on the leaves
$\{w_2=c\}$, for all sufficiently small
$c\in\bC$. Again, this result is optimal, because
the biholomorphism 
$$
H(z_1,z_2,w_1,w_2)=(z_1e^{iw_2},z_2,w_1,w_2)
$$  maps the set $A$ into itself. This map is only
algebraic on the leaves $\{w_2=c\}$ and not in
the whole space.
\medskip

\noindent{\bf Example 3.1.7.} Consider the submanifold
of $\bC^4$ defined by
$$
\left\{\aligned \re w_1 &=|z_1|^2\\\im w_1 &=
|z_2|^2\\ \im w_2 &=0.
\endaligned\right.\tag3.1.6
$$ Note that this submanifold is not generic (nor is
it CR!) on the set
$\{(z_1,z_2,w_1,w_2)=(0,0,0,s_2)\:s_2\in\bR\}$.
However, outside that set the manifold
\thetag{3.1.6} is a generic holomorphically
nondegenerate manifold $M$. The function
$h(z,w)=w_2$ is real on $M$, and
$M\cap\{w_2=s_2\}$ is generically minimal for all
$x\in
\bR$. As above, Theorem 3.1.2 implies that any
holomorphic mapping
$H\:\bC^4\mapsto
\bC^{N^\prime}$ near 0, generically of rank 4, such
that $H(M)$ is contained in a 5 dimensional real
algebraic subset of
$\bC^{N^\prime}$ is algebraic on the leaves
$\{w_2=c\}$, for all sufficiently small
$c\in\bC$. We invite the reader to construct an
example, e.g. similar to the ones considered above,
to show that one cannot have a stronger conclusion.
\medskip

We can also formulate a result 
 that holds at most, but
not necessarily all, points of the algebraic
set.
\proclaim{Theorem 3.1.8} Let $A\subset \cnn$ be
an irreducible, holomorphically nondegenerate,
real algebraic set, and let $\scrV$ be a
complex algebraic variety in
$\bC^N$ that contains $A$. Then either of the
following holds, for all points $p\in \areg$
outside a proper real algebraic subset of
$A$: 
\roster
\item"(i)" There is
$h\in\scrA_\scrV(p)$ such that $h$ is not constant
and $h|_A$ is real valued.
\item"(ii)" All mappings
$H\in\hol_\scrV(p,p^\prime)$, where
$p^\prime\in\bC^ {N^\prime}$ is arbitrary, such that
the generic rank of $H$ equals
$\dim_\bC\scrV$ and such that
$H(A)$ is contained in a real algebraic set
$A^\prime$, with $p^\prime\in A^\prime$ and
$\dim_\bR A=\dim_\bR A^\prime$,
 are algebraic in $\scrV$, i.e.
$H\in\scrA_\scrV(p,p^\prime)$.
\endroster
\endproclaim

Before we proceed with the proofs of Theorems 3.1.2
and 3.1.8 (\S3.3 and
\S3.4), we need a result on ``propagation of
algebraicity'' that we establish in the next
subsection. 

\subhead 3.2. Propagation of algebraicity\endsubhead
We assume that we have an algebraic foliation of  some
domain in complex space, and a holomorphic function 
$f$ whose restriction to a certain
sufficiently large collection of the leaves is
algebraic. We shall show that the restrictions of $f$
to all leaves in the domain are also algebraic,
provided that the domain has a nice ``product
structure'' with respect to the foliation. This will be
essential in the proof of Theorem 3.1.2. This result
may already be known.
\proclaim{Lemma 3.2.1} Let $f(z,w)$ be a holomorphic
function in
$U\times V$, where $U\subset\bC_z^a$ and
$V\subset\bC_w^b$ are domains. Assume that there is
a subdomain
$V_0\subset V$ and a nontrivial polynomial 
$P(z,X;w)\in\scrO_b(V_0)[z,X]$, i.e. 
$P$ is a polynomial in
$z=(z_1,...,z_a)$ and $x$ with coefficients
holomorphic in $V_0$, such that
$$ P(z,f(z,w);w)\equiv0, \ \ z\in U,\ w\in
V_0. \tag3.2.1
$$ Then there is a nontrivial
polynomial 
$\tilde P(z,x;w)\in\scrO_b(V)$
$[z,X]$,  such that
$$
\tilde P(z,f(z,w);w)\equiv0,  \ \ z\in U,\ w\in
V.\tag3.2.2
$$
\endproclaim
\demo{Proof} Pick any point $w^0\in V_0$, and
consider $P=P(z,X;w)$ as an element of
$\scrO_b(w^0)[z,X]$. We order the monomials
$z^\alpha$ by choosing a
bijection $i:\bZ^a_+\to \bZ_+$, and  write 
$$ P(z,X;w)=\sum_{k=0}^pp_k(z;w)X^k,\tag3.2.3
$$ where each
$p_k(z;w)\in\scrO_b(w^0)[z]$ is of the form
$$ p_k(z;w)=\sum_{i(\alpha)\leq
q_k}a^k_\alpha(w)z^\alpha\tag3.2.4
$$ with
$a^k_\alpha\in\scrO_b(w^0)$. We choose
$p,q_1,...,q_p$ minimal such that $P$ can be written
in this form with the leading terms in
\thetag{3.2.3} and
\thetag{3.2.4}  not identically 0.  We may
assume that the numbers $p$ and $q_p$ are minimal in
the sense that if
$P^\prime$ is another polynomial in
$\scrO_b(w^0)[z,x]$, with corresponding numbers $p
^\prime$ and $q^\prime_{p^\prime}$, such that
\thetag{3.2.1} holds then
$p\leq p^\prime$ and if
$p=p^\prime$ then $q_p\leq q^\prime_{p^\prime}$. The
polynomial $P$ is then unique  modulo multiplication
by elements in $\scrO_b(w^0)$ in the sense that
if $P^\prime$ is as above with
$p^\prime=p$ and
$q^\prime_{p^\prime}=q_p$ then there are germs
$c_1(w),c_2(w)\in\scrO_b(w^0)$, not identically 0,
such that
$$ c_1(w)P(z,x;w)\equiv
c_2(w)P^\prime(z,x;w).\tag3.2.5
$$ Since $p$ is minimal, the function $p_0(z;w)$ is
not identically 0, and thus  there is a coefficient
$a^0_{\alpha_0}(w)$ which is not identically 0.
 The equation \thetag{3.2.1} can then be written
in the form
$$ Q(z,f(z,w);w)\equiv
-a^0_{\alpha_0}(w)z^{\alpha_0},\tag3.2.6
$$
with $Q(z,X;w)\in \scrO_b(w^0)[z,X]$. Now, the
uniqueness of $P$ in the sense of
\thetag{3.2.5} and the fact that
$a^0_{\alpha_0}\not\equiv0$ imply that the
coefficients of $Q(z,X;w)$ satisfying (3.2.6) are
actually unique. After dividing (3.2.6) by
$-a^0_{\alpha_0}(w)$, we find $Q'(z,X;w)\in
\scrM_b(w^0)[z,X]$ satisfying
$$
Q'(z,f(z,w);w)\equiv
z^{\alpha_0},\tag3.2.7
$$
where $\scrM_b(w^0)$ denotes the field of
meromorphic functions near $w_0$. 

 We order the
set of indices
$(k,\alpha)$, for $k=0,...,p$ and
$i(\alpha)\leq q_k$, minus the index $(0,\alpha_0)$
in some way, e.g.  the ``canonical'' way induced by
the ordering $i=i(\alpha)$. We hence obtain a
bijection
$(k,\a) \mapsto j(k,\alpha)$ from this set of
indices to the set of numbers
$\{1,2,...,\mu\}$, where $\mu$ is the number of
elements in this set of indices. We introduce the
vector valued functions
$A(z;w)\in\bC^\mu$ defined by letting the $j$th
component be
$$ A_{j}(z;w)=z^\alpha f(z,w)^k,
\ \ \text {for}  \ j = j(k,\alpha),\tag3.2.8
$$ and $b(w)\in\bC^\mu$ defined by
$$
b_{j}(w)={a^k_{\alpha}(w)\over
-a^0_{\alpha_0}(w) }, \ \ \text {for}  \ j =
j(k,\alpha) .\tag3.2.9
$$ Then \thetag{3.2.7} can be written
$$ A(z;w)\cdot b(w)\equiv
z^{\alpha_0},\tag3.2.10$$ where $\cdot$ denotes the
usual dot product of vectors in $\bC^\mu$.
Moreover, the vector valued meromorphic function
$b(w)$ is the  unique meromorphic solution of
\thetag{3.2.10}.  Consider the
$\mu\times\mu$ matrix valued holomorphic function
$B(z_1,...,z_\mu,w)$ (of $a\mu+b$ variables) defined
by letting the matrix element
$B_{ij}(z_1,...,z_\mu;w)$, for
$i,j=1,...,\mu$, be
$$ B_{ij}(z_1,...,z_\mu;w)=A_j(z_i;w).\tag3.2.11$$
We claim that the determinant of
$B(z_1,...,z_\mu;w)$ is not  identically 0.
Indeed, if it were, then we could find a vector
valued holomorphic function $c(w)$, not
identically 0, such that $ A(z;w)\cdot
c(w)\equiv0$, which would contradict the uniqueness
of the solution $b(w)$ of (3.2.10).  Thus, we can
find fixed values 
$z_1^0,...,z_\mu^0$ such that
$\Delta(w)$, the determinant of
$B(z_1^0,..., z_\mu^0;w)$ as a function of $w$, is
not identically 0. We can then solve for $b(w)$ as
the unique solution of the system obtained from
\thetag{3.2.10} after substituting successively
$z_1^0,...,z_\mu^0$ for $z$. Since the matrix
$B(z_1^0,..., z_\mu^0;w)$ has entries holomorphic
in all of $V$, by Cramer's rule it follows that
the solution 
$b(w)$ thus obtained is in
$\scrM_b(V)$.  
 Hence  $Q'(z,X;w) \in
\scrM_b(V)[z,X]$. After clearing
denominators we obtain (3.2.2) from (3.2.7). This
completes the proof of Lemma 3.2.1.
\qed
\enddemo

\subhead 3.3. Proof of Theorem 3.1.2\endsubhead First, since
all assumptions and conclusions in the theorem are
related to $\scrV$, and
$p_0\in\scrV$ is a regular point of
$\scrV$, it suffices to consider the
case where $\scrV=\bC^N$ and
$M$ is generic; we will assume this for the
rest of the proof.  By assumption (i) in the theorem, we
can find algebraic coordinates
$(u,v)\in\bC^{N-q}\times\bC^q$,  vanishing at $p_0$,
in a neighborhood $U_1$ of $p_0$ such that $h_j=v_j$
for $j=1,...,q$. We may assume
$U_1=A_1\times
B_1$ , where
$A_1\subset\bC^{N-q}$ and $B_1
\subset\bC^q$. It follows from the assumptions that
 $M\cap S_{h(p)}$
is minimal at $p$ for $p$  outside a proper real
algebraic subset of
$M\cap U_1$. Similarly,
$M$ is
$l(M)$-nondegenerate outside a proper real algebraic
subset of
$M$, where
$l(M)$ is the Levi number defined in \S1.3. 
Also, the mapping $H$ attains its maximal
rank outside a proper complex analytic subset of $\bC^N$
near $p_0$, and  since $M$ is  generic it is not
contained in any proper complex analytic set. Thus,
$H$ attains its maximal rank at points on $M$ outside
a proper real  analytic subset of $M$. Finally, for each
$j$, the $j$th Segre set $N_j(Z)$  of
$M$ at
$Z$ (defined in \S 2.2) has maximal generic dimension
for
$Z$ outside a proper real algebraic subset of $M$. Hence
we can find
$p_1 \in M \cap U_1$ such that
\roster
 \item "(a)" $M \cap S_{h(p_1)}$ is minimal at $p_1$,
\item "(b)" $H$ has rank $N$ at $p_1$,
\item "(c)" $M$ is $l(M)$-nondegenerate at $p_1$,
\item "(d)" For each $j$, the generic dimension $d_j$ of 
$N_j(p_1)$ is maximal.
\endroster 
 We will prove Theorem 3.1.2
by first showing that there is a neighborhood of
$p_1$ in
$\bC^N$ such that $H|_{S_{h(p)}}$ is algebraic for  every
$p$ in that neighborhood, and then applying Lemma 3.2.1
to deduce the full statement of the theorem.
For this, we claim that we may assume that
the target
$A^\prime$ is contained in $\bC^N$  and that $H$ is a
mapping into $\bC^N$. Indeed, there is a neighborhood
$U_2
\subset U_1$ of the point $p_1$ such that
$Y^\prime=H(U_2)$ is a complex  holomorphic
submanifold of dimension $N$ in $\bC^{N^\prime}$
through the point 
$p_1^\prime=H(p_1)$ . Since
$M$ is generic and $H$ is a biholomorphism of $U_2$
onto
$Y^\prime$, it follows that
$A^\prime$ is a generic submanifold of $Y^\prime$
near $p_1^\prime$ .  Denote by $M^\prime$ a piece of
$A'$ near $p_1^\prime$ and choose it such that $M'$ is a
generic submanifold of $Y'$. Then,
$M^\prime$ is real algebraic and its intrinsic
complexification 
$\scrV^\prime\subset\bC^{N^\prime}$ is a complex
algebraic manifold near $p_1^\prime$. Since both
$Y^\prime$ and
$\scrV^\prime$ contain $M^\prime$, and $M^\prime$ is
generic in both manifolds, it follows that $Y^\prime
=
\scrV^\prime$.  We can therefore choose algebraic
coordinates in a neighborhood
$U^\prime_2$ of $p_1^\prime$ in
$\bC^{N^\prime}$, vanishing at
$p_1^\prime$, such that $H=(\hat H,0)$ in these
coordinates and
$\hat H$ maps $M\cap  U_2$ into
$M^\prime\cap U_2^\prime\subset Y^\prime\cong\bC^N$. 
 In what follows, we assume
that 
$\bC^{N^\prime}=\bC^N$ and we take
$\hat H$ as our mapping $H$. 

Let
$(z,w)\in\bC^n\times\bC^d=\bC^N$, where $n$ is the CR
dimension and
$d$  the codimension of $M$, be (algebraic) normal
coordinates for
$M$, vanishing  at $p_1$, i.e. $M$ is defined near
$p_1$ by
\thetag{1.1.3} and similarly for the target
$M^\prime$ (denoting the function defining $M^\prime$
by
$Q^\prime$). We write
$(z,w)=(z(u,v),w(u,v))$ to indicate the relationship
between the local normal coordinates $(z,w)$ near
$p_1$ and the coordinates $(u,v)$ in $U_1$. Thus, we
can write the mapping $H$ as 
$H=(f,g)$, where $f(z,w)\in\bC^n$ and
$g(z,w)\in\bC^d$, such that
$$
\bar g=\bar Q^\prime(\bar f,f,g)\tag3.3.1
$$ holds for points $(z,w)\in M$ near $p_1=0$.  By
complexifying, we obtain
$$
\bar g (\chi,\tau)=\bar Q^\prime(\bar f
(\chi,\tau),f(z,w),g(z,w)),\tag3.3.2
$$ for all
$(z,w,\chi,\tau)\in\scrM$ near 0.  We
define the holomorphic vector fields  $\Cal L_j$ in
$\bC^{2N}$ tangent to  $\scrM$
(and corresponding to the CR vector fields of $M$) by
$$
\Cal L_j=\frac{\partial}{\partial
\chi_j}+\sum_{k=1}^d \bar Q_{k,\chi_j}(\chi,
z,w)\frac{\partial}{\partial
\tau_k}, \ \ j = 1,\ldots, n.\tag3.3.3
$$ 
 We
shall also need 
the following vector fields  tangent to
$\scrM$ given by
$$\aligned \tilde\Cal
L_j&=\frac{\partial}{\partial
z_j}+\sum_{k=1}^d  Q_{k,z_j}(z,
\chi,\tau)\frac{\partial}{\partial
w_k}, \quad\quad j=1,...,n,\\ \Cal
T_j&=\frac{\partial }{\partial w_j}+\sum_{k=1}^d
\bar Q_{k,w_j}(\chi,
z,w)
\frac{\partial}{\partial
\tau_k},\quad\quad j=1,...,d,\\
V_j&=\tilde\scrL_j-\sum_{k=1}^d Q_{k,z_j}(z,\chi,
\tau)\Cal T_k,  \quad\quad j=1,...,n.
\endaligned\tag3.3.4
$$  Note that the coefficients of all
the vector fields given by (3.3.3) and (3.3.4) are
algebraic functions of
$(z,w,\chi,\tau)$.
\proclaim{Assertion 3.3.1} There is a neighborhood
$U_3\subset U_2$ such that, for all
$(z,w,\chi,\tau)\in\scrM\cap(U_3\times {}^*U_3)$ and
all multi-indices
$\gamma=(\gamma^\prime,\gamma^{\prime\prime})$, 
$$
\frac{\partial^{|\gamma|}f_j}{\partial
z^{\gamma^\prime}\partial w^{\gamma^ {\prime\prime}}}
(z,w)=\Psi^\gamma_j(...,V^{\alpha^3}\Cal T^{\alpha^2}
\Cal L^{\alpha^1}\bar f_k(\chi,
\tau),..., V^{\beta^3}\Cal T^{\beta^2}\Cal L^{\beta^1}
\bar g_l(\chi,\tau),...),\tag3.3.5
$$ where $j,k=1,...,n$,
$l=1,...,d$,
$|\alpha^1|,|\beta^1|\leq l(M)$,
$|\alpha^2|,|\beta^2|\leq|\gamma^{\prime\prime}|$, 
$|\alpha^3|,|\beta^3|\leq|\gamma^\prime|$, and the 
$\Psi^\gamma_j$ are algebraic holomorphic functions
of their  arguments.\endproclaim
\demo{Proof} We apply the operators
$\Cal L_j$ to the identity
\thetag{3.3.2}, and use the fact that the matrix
$\Cal L\bar f$ at
$(z,w,\chi,\tau)=(0,0,0,0)$ is invertible (since $H$
is a biholomorphism at $p_1=0$) to deduce that there
are algebraic functions $F_j$ such that, for points
on 
$\scrM$ near 0,
$$ \bar Q^\prime_{\chi_j}(\bar f,f,g)=F_j(\Cal L\bar
f,\Cal L\bar g).\tag3.3.6
$$ 
 We repeat this procedure, using in the
next step \thetag{3.3.6} instead of 
\thetag{3.3.2} and so on. Since
$H$ is a  biholomorphism at $p_1$,
$M^\prime$ is $l(M)$-nondegenerate at $p_1^\prime$.  Hence
(see \S 1.3)
$$
\text{span}\{\bar
Q^\prime_{z,\chi^{\alpha}}(0,0,0)\:|\alpha|\leq
l(M)\}=\bC^n.
\tag3.3.7$$
It follows from the algebraic implicit function
theorem  and (3.3.7)
 that, for all
$(z,w,\chi,\tau)\in\scrM$ near the origin,
$$ f_j(z,w)=\Psi_j(...,\Cal L^{\alpha}\bar
f_k(\chi,\tau),...,\Cal L^{\beta}
\bar  g_l(\chi,\tau),...)\quad,\quad
j=1,...,n,\tag3.3.8
$$ where $k=1,...,n$, $l=1,...,d$,
$|\alpha|,|\beta|\leq l(M)$, and the 
$\Psi_j$ are algebraic holomorphic functions of
their  arguments (cf. e.g. [BR4, Lemma 2.3]). Now,
since
$f(z,w)$ is a function of
$(z,w)$ only, we have, for any multi-index
$\gamma=(\gamma^\prime,\gamma^ {\prime\prime})$,
$$ V^{\gamma^\prime}\Cal
T^{\gamma^{\prime\prime}}f(z,w)=\frac{\partial^
{|\gamma|}f_j}{\partial z^{\gamma^\prime}\partial
w^{\gamma^{\prime\prime}}} (z,w).\tag3.3.9
$$ The assertion follows if we apply
$V^{\gamma^\prime}\Cal T^{\gamma^{\prime\prime}}$ to
\thetag{3.3.8}, since the $V_j$ and
$\Cal T_l$ are tangent to $\scrM$. \qed\enddemo 

We
now proceed with the proof of Theorem 3.1.2. From
(3.3.2) we have 
$$ g_l(z,w)=Q^\prime_l(f(z,w),\bar f(\chi,\tau),\bar
g(\chi,\tau))\tag3.3.10
$$  for $(z,w,\chi,\tau)\in\scrM$ and $l=1,...,d$. If
we apply 
$V^{\gamma^\prime}\Cal T ^{\gamma^{\prime\prime}}$ to
this equation we obtain
$$
\aligned &\frac{\partial^{|\gamma|}g_l}{\partial
z^{\gamma^\prime}\partial w^{\gamma^ {\prime\prime}}}
(z,w)=\\&\hskip .3truecm\Phi^\gamma_l\left(...,
\frac{\partial^{|\alpha^1|+|\alpha^2|}f}{\partial
z^{\alpha^2}\partial w^ {\alpha^1}}(z,w),...,
V^{\beta^2}\Cal T^{\beta^1}
\bar f_j(\chi,
\tau) ,..., V^{\mu^2}\Cal T^{\mu^1}
\bar g_k(\chi,\tau),...\right),\endaligned\tag3.3.11
$$ where $j=1,...,n$,
$k,l=1,...,d$,
$|\alpha^1|,|\beta^1|,|\mu^1|\leq
|\gamma^{\prime\prime}|$,
$|\alpha^2|,|\beta^2|,|\mu^2|\leq|\gamma^{\prime}|$,
and where $\Phi^\gamma  _l$ are algebraic holomorphic
functions of their arguments. Using \thetag {3.3.5},
we obtain
$$
\frac{\partial^{|\gamma|}g_l}{\partial
z^{\gamma^\prime}\partial w^{\gamma^ {\prime\prime}}}
(z,w)=\Xi^\gamma_l(...,V^{\alpha^3}\Cal T^{\alpha^2}
\Cal L^{\alpha^1}\bar f_j(\chi,
\tau),..., V^{\beta^3}\Cal T^{\beta^2}\Cal L^{\beta^1}
\bar g_k(\chi,\tau),...),\tag3.3.12
$$ where $j=1,...,n$,
$k,l=1,...,d$,
$|\alpha^1|,|\beta^1|,|\mu^1|\leq l(M)$,
$|\alpha^2|,|\beta^2|,|\mu^2|\leq
|\gamma^{\prime\prime}|$,
$|\alpha^3|, |\beta^3|,|\mu^3|\leq |\gamma^\prime|$,
and the
$\Xi^\gamma_l$ are holomorphic algebraic functions of
their arguments. For notational brevity, we use the
notation $Z=(z,w)$ and
$\zeta=(\chi,\tau)$. If we denote by
$H(Z)=(H_1(Z),..., H_N(Z))$ the components of $H$ in
an arbitrary algebraic coordinate system near the
point $p_1^\prime=H(p_1)$  then it follows
from (3.3.5) and \thetag{3.3.12} that we have
$$
\frac{\partial^{|\gamma|}H_k}{\partial
Z^\gamma}(Z)=\Theta^\gamma_k
(Z,\z,\ldots, \frac{\partial^{\alpha}\overline H_j}{\partial
\z^\alpha}(\z),\ldots),
\tag3.3.13
$$ where $j,k=1,...,N$, $\gamma$ arbitrary, $|\alpha| \le
|\gamma| + l(M)$, and
$\Theta^\gamma_k$ are holomorphic algebraic functions
 of
their arguments for $(Z,\z) \in \scrM$ near $(p_1,\bar
  p_1)$.
\proclaim{Assertion 3.3.2} For $Z \in M$ near $p_1$,
 let $N_j(Z)$ denote the $j$th Segre set of $M$ at
$Z$ and $d_j$ the generic dimension of $N_j(p_1)$. For
some $j$, $1 \le j\le j_0-1$, let 
$$
\bC^{d_j}\times \bC^N\ni (s, Z)\mapsto\zeta
(s,Z)\in\bC^N,\tag3.3.14
$$ be an algebraic map,
holomorphic near
$(0,p_1)$.  Suppose
 $s\mapsto
\zeta(s,p_1)$ has generic rank $d_j$, and  $\zeta(s,
Z)\in{}^*N_j(Z)$. Then there is an algebraic map
$$
\bC^{d_{j+1}}\times\bC^N\ni
(t,Z)\mapsto(\Pi(t,Z),s(t))\in
\bC^N\times\bC^ {d_j},\tag3.3.15
$$ holomorphic near $(0,p_1)$, such that 
$ (\Pi(t,Z),\zeta(s(t),Z))\in\scrM$, 
the mapping $t\mapsto \Pi(t,Z)$ has generic rank
$d_{j+1}$, and $\Pi(t,Z)\in N_{j+1}(Z)$, for all
$Z\in M$ near $p_1$.\endproclaim
\noindent{\bf Remark:} If
$\zeta(s,Z)$ is algebraic anti-holomorphic in $Z$
rather than algebraic holomorphic then the same
conclusion holds with ``holomorphic" replaced by
``anti-holomorphic".\medskip

\demo{Proof} Note first that by assumption (d), $d_j$
is also the generic dimension of $N_j(Z)$ for $Z$ near
$p_1$.  We write the map
$\zeta(s,Z)$ in the normal coordinates as
$(\chi(s,Z),\tau(s,Z))$. For $Z \in M$ fixed near
$p_1$ consider the map
$$
\bC^n\times\bC^{d_j}\ni(z,s)\mapsto(z,Q(z,\chi(s,Z),
\tau(s,Z)))
\in\bC^n\times\bC^d=\bC^N.\tag3.3.16
$$ Note that
$
(z,Q(z,\chi(s,Z),\tau(s,Z)),\chi(s,Z),\tau(s,Z))\in\scrM
$. Since $N_{j+1}(Z)$, for
$Z\in M$, is defined as
$
\{(z,Q(z,\chi,\tau))\:\exists(\chi,\tau)\in{}^*N_j(Z)\}
$, and the mapping
$s\mapsto\zeta(s,Z)\in{}^*N_j(Z)$ has  rank $d_j$,
which is also the generic dimension of $^*N_j(Z)$, it is
easy to verify that the map
\thetag{3.3.16}  has generic rank
$d_{j+1}$. Thus, by the rank theorem, there is
an algebraic map 
$$
\bC^{d_{j+1}-n}\ni t^\prime\mapsto
s(t^\prime)\in\bC^{d_j}\tag3.3.17
$$ such that
$$
\bC^n\times\bC^{d_{j+1}-n}\ni(z,t^\prime)=t\mapsto(z,Q(z,\chi(s(t^\prime),
Z),\tau(s(t^\prime),Z)))\in\bC^N\tag3.3.18
$$ has rank $d_{j+1}$. The proof of Assertion 3.3.2
follows by taking
$t=(z,t^\prime)$, 
$s(t)=s(t^\prime)$, and,
$$
\Pi(t,Z)=(z,Q(z,\chi(s(t^\prime),Z),
\tau(s(t^\prime),Z))).\tag3.3.19
$$
\qed\enddemo Now, define the map
$\Pi_0(Z)=Z$ and the map
$\zeta_0(Z)=\bar Z$.  The latter, thought of as a
map
$\bC^0\times\bC^N\mapsto \bC^N$, satisfies the
hypothesis of Assertion 3.3.2 above with $j=0$ (see
the remark following the assertion). Thus, we get an
algebraic map, anti-holomorphic in $Z$
$$
\bC^{d_1}\times\bC^N\ni (t,Z)\mapsto
\Pi_1(t,Z)\in\bC^N,\tag3.3.20
$$ of rank $d_1$ in $t$ for $Z$ near $p_1\in M$,
such that
$ (\Pi_1(t,Z),\bar Z)\in\scrM
$, and  $\Pi_1(t,Z)\in N_1(Z)$ for
$Z\in M$. Defining the map $\zeta_1(t,Z)$ by 
$$
\zeta_1(t,Z)=\overline{\Pi_1(\bar t,Z)}\tag3.3.21
$$ we obtain a map into
$^*N_1(Z)$ that satisfies the hypothesis of the
assertion with
$j=1$  (this time the map is holomorphic in $Z$).
Applying the assertion again and proceeding
inductively, we obtain a sequence of algebraic  maps 
$$ \Pi_0(Z),\Pi_1(t_1,Z),...,\Pi_{j_0}(t_{j_0},Z),
\zeta_0(Z),
\zeta_1(s_1,Z),...,\zeta_{j_0-1}(s_{j_0-1},Z)
$$  (either holomorphic or anti-holomorphic in $Z$)
with
$t_j\in\bC^{d_j}$,
$s_j\in\bC^{d_j}$, and accompanying maps
$s_1(t_2),...,s_{j_0-1}(t_{j_0})$ such that the maps
$ t_j\mapsto \Pi_j(t_j,Z)$ and $s_j\mapsto
\zeta_j(s_j,Z)
$ are of rank $d_j$, map into
$N_j(Z)$ and $^*N_j(Z)$ respectively for $Z\in
M$, and satisfy
$$
(\Pi_{j+1}(t_{j+1},Z),\zeta_j(s_j(t_{j+1}),Z))
\in\scrM,\tag3.3.22
$$ for $j=0,...,j_0-1$. Morever, we have the relation
$$
\zeta_j(s_j,Z)=\overline{\Pi_j(\bar
s_j,Z)}.\tag3.3.23
$$
\proclaim{Assertion 3.3.3} For each
$j=1,...,j_0$, 
$$
\frac{\partial^{|\gamma|}H_k}{\partial
Z^\gamma}(\Pi_j(t_j,Z))=F^\gamma_{jk}
\left(t_j,Z,\bar
Z,...,\frac{\partial^{|\alpha|}H_l}{\partial
Z^\alpha} (Z), ...,\frac{\partial^{|\beta|}\bar
H_l}{\partial \zeta^\beta}(\bar Z),...
\right)
\tag3.3.24
$$ holds for $Z\in M$ near $p_1$,  where
$k,l=1,...,N$,
$|\alpha|,|\beta|\leq |\gamma|+j\,l(M)$, and
$F^\gamma_{jk}$ are holomorphic algebraic functions of
their arguments. 
\endproclaim
\demo{Proof} The proof is by induction on $j$. For
$j=1$, we prove the statement by taking
$Z$ to be $\Pi_1(t_1,Z)$ and
$\zeta$ to be $\z_0(Z)=\bar Z$ in
\thetag{3.3.13} (using
\thetag{3.3.22}). Assume now that (3.3.24) holds for
$j=1,...,i$ (with $i<j_0$). By
\thetag{3.3.23} we have
$$
\frac{\partial^{|\gamma|}\bar H_k}{\partial
\zeta^\gamma}(\zeta_i(s_i,Z))=
\bar F^\gamma_{jk}
\left(s_j,\bar
Z,Z,...,\frac{\partial^{|\alpha|}\bar
H_l}{\partial \zeta ^\alpha}(\bar Z),
...,\frac{\partial^{|\beta|} H_l}{\partial
Z^\beta}(Z),...\right).
\tag3.3.25
$$ Now (3.3.24) follows for $j=i+1$ from (3.3.25)
by taking
$Z$ to be $\Pi_{i+1}(t_{i+1},Z)$ and
$\zeta$ to be $\zeta_i(s_i(t_{i+1}),Z)$ in
\thetag{3.3.13}.  \qed\enddemo

We now complete the proof of Theorem 3.1.2. For $p$ near
$p_1$ it follows from  Corollary 2.2.2, since 
$M\cap S_{h(p)}$ is minimal, that the
maximal Segre set $N_{j_0}(p)$ is contained in and
contains an open piece of
$S_{h(p)}$. Since $M$ is generic, it is easy to see
that $h(M)$ contains an open neighborhood  of
$c^1=h(p_1)$ in $\bR^q$. Thus, by the rank theorem, and
using the coordinates $(u,v)$ in $U_1$, there is a real
algebraic injective map
$\bR^q\ni c\mapsto(u(c),c)\in M$, for $c$ near $c^1$, which
can be complexified to  an algebraic injective map
$v\mapsto (u(v),v)$, for $v$ in a neighborhood of
$c^1$ in $\bC^q$. Now, let $Z$ be the  point
$Z(c)=(z(u(c),c), w(u(c),c))$ where $c\in\bR^q$ is
some arbitrary point near $c^1$. Applying  Assertion
3.3.3 with this choice of
$Z$, $\gamma=0$ and $j=j_0$, we deduce that each
component $H_l$ is algebraic on $S_c$ and  satisfies
there a polynomial equation with coefficients that
depend real analytically on $c$ (we may take $t_{j_0}$
as algebraic coordinates on $S_c$). In terms of the
coordinates
$(u,v)$ with
$\tilde H(u,v)=H(Z(u,v))$, there are polynomials
$P_l(u,X;c)$ in
$(u,X)\in\bC^{N-q}\times\bC$,
$l=1,...,N$,  with coefficients that  are real
analytic functions in $c$, for $c$ close to $c^1$,
such that
$$ P_l(u,\tilde H_l(u,c);c)\equiv0\tag3.3.26
$$ ($u$ are also algebraic coordinates on $S_c$ and
it is easy to see that the algebraic change of
coordinates $u=u(t_{j_0})$ on $S_c$ depends real
algebraically on $c$). Extending the coefficients of
the polynomials to be complex analytic functions of
$v$ in a neighborhood
$B_2$ of $c^1$ in $\bC^q$, we obtain  polynomials
$P_l(u,X;v)\in
\scrO_q(B_2) [u,X]$ such that 
$$ P_l(u,\tilde H_l(u,v);v)\equiv0,\tag3.3.27
$$ holds in $A_1\times B_2$. Since
$\tilde H_l(u,v)$ is holomorphic in 
$A_1\times B_1$, there is, by Lemma 3.2.1, possibly
another polynomial
$\tilde P_l(u,X;v)
\in\scrO_q(B_1)[u,X]$ such that
$$
\tilde P_l(u,\tilde H_l(u,v);v)\equiv0,\tag3.3.28
$$ holds in $U_1=A_1\times B_1$. This completes the
proof of Theorem 3.1.2 with $U=U_1$, and $\delta>0$ 
being any number such that the ball of radius
$\delta$ centered at
$v=0$ is contained in $B_1$ (recall that $B_1$  is a
neighborhood of
$v=0$). The proof of Theorem 3.1.2 is now complete.\qed

\subhead 3.4. Proof of Theorem 3.1.8\endsubhead
Put $M=\areg$. First, note that if $M$ is contained
in a proper complex algebraic subvariety of $\scrV$
then (i) holds for all points $p\in M$. If 
$M$ is not contained in a proper algebraic subvariety
of $\scrV$ then
$M$ is a generic real algebraic submanifold of
$\scrV_{\text{reg}}$ at $p$, for all $p$ outside some
proper real algebraic subvariety of $M$. Thus, as in
the proof of Theorem 3.1.2, we may assume that
$\scrV=\bC^N$ and that $M$ is a generic
holomorphically nondegenerate submanifold in
$\bC^N$.  Let $\po \in M$ be a point whose
 CR orbit
 has maximal dimension. If $M$ is
minimal at $p_0$ then  (ii) holds with $p=p_0$,  by
Corollary 3.1.4.  Moreover, if $M$ is minimal at
$\po$
 then $M$ is minimal for $p$ outside a real algebraic
variety and therefore  (ii) holds at such $p$. Thus, the
theorem follows if we can show that
$M$ is minimal at $p_0$ unless (i) holds at $p_0$ . 
The proof of Theorem 3.1.8 will then be
completed by the following lemma.
\proclaim{Lemma 3.4.1} Let $M$ be a generic real
algebraic submanifold in $\bC^N$, and let $p_0\in M$
with CR orbit
of maximal dimension. Then $M$ is minimal at 
$p_0$ if and only if there is no, non-constant 
$h\in\scrA_N(p_0)$ such that $h|_M$ is real valued.
More  precisely, if the codimension of the local CR
orbit of $p_0$ in
$M$ is $q$ then there are
$h_1,...,h_q\in \scrA_N(p_0)$ such that $h_j|_M$ is 
real valued for
$j=1,...,q$ and
$$
\partial h_1(\po)\wedge...\wedge\partial
h_q(\po)\neq0.\tag3.4.1
$$\endproclaim
\noindent{\bf Remark:} Lemma 3.4.1 implies that the
decomposition of $M$ into CR orbits near $p_0$
 is actually an algebraic foliation,
because the CR orbit of a point
$p_1$ near $p_0$ must equal $\{p\in
M\:h_j(p)=h_j(p_1)\,,\, j=1,...,q\}$.  If Corollary
2.2.5 is viewed as an algebraic version of the Nagano
theorem (for the special class of algebraic vector
fields that arise in this situation; see the paragraph
following Corollary 2.2.5) then this lemma is the
algebraic version of the Frobenius theorem. \medskip
\demo{Proof of Lemma {\rm 3.4.1}} Assume that there
is a non-constant
$h\in\scrA_N(p_0)$ such that
$h|_M$ is real. Then, by Lemma 3.1.1,
$M\cap\{Z\:h(Z)=h(p)\}$ is  a CR submanifold for all
$p\in M$ near $p_0$ such that $dh(p)\neq0$. Since
$h$ is real on $M$ all CR and  anti-CR vector fields tangent
to
$M$ annihilate $h$; hence  the submanifold 
$M\cap\{Z\:h(Z)=h(p)\}$ has the same CR dimension as
$M$.
Thus,
$M$ is not minimal at $p$. Since this is true for all
$p\in M$ near
$p_0$ outside a proper real algebraic subset, 
$M$ is not minimal anywhere. This proves the ``only
if''  part of the first statement of the lemma. The
``if'' part will follow from the more precise
statement at the end of the lemma, which we shall
now prove.

We choose algebraic normal coordinates
$(z,w)\in\bC^N$,  vanishing at
$p_0$, such that $M$ is given by \thetag{1.1.3} near
$p_0$. Denote by $W_0$ the  CR orbit of
$p_0=0$, and by
$X_0$ its intrinsic complexification. By Theorem
2.2.1, $N_{j_0}(\po)$, the maximal Segre set of 
$M$ at $p_0$,  is contained in and contains
an open piece of $X_0$. The complex dimension of
$X_0$ is
$d_{j_0}$, the generic dimension of $N_{j_0}(\po)$. Since
the codimension of
$W_0$ in 
$M$ is
$q$, the complex codimension of its intrinsic
complexification $X_0$ is also $q$, i.e.
$d_{j_0}=n+d-q$. Let $r=d-q$. By a
linear change of the $w$ coordinates, we may assume
that the tangent plane of
$X_0$ at $0$ is
$\{(z,w)\:w_{r+1}=...=w_d=0\}$. We decompose
$w$ as
$(w^\prime,w^{\prime\prime})\in\bC^r\times\bC^{q}=\bC^d$.
Note  that at the  point $\tilde p=(0,s)\in M$, where
$s=(s^\prime,s^{\prime\prime})\in\bR^r
\times\bR^{q}$, 
$(\tilde z,\tilde w)=(z,w-s)$ are normal coordinates
vanishing at
$\tilde p$ and $M$ is given by
$$
\tilde w=Q(\tilde z,\bar{\tilde z},\bar{\tilde
w}+s)-s.\tag3.4.2
$$ We denote by
$W_{s^{\prime\prime}}$ the local CR orbit of $(0,0,s^
{\prime\prime})$, by $X_{s^{\prime\prime}}$ its
intrinsic complexification, and  by
$N_{j_0}(s^{\prime\prime})$ the maximal Segre set at $(0,0,s^
{\prime\prime})$ .
Since the CR orbit at
$p_0$ has maximal dimension, all
$W_{s^{\prime\prime}}$,
$X_{s^{\prime\prime}}$, and
$N_{j_0} (s^{\prime\prime})$ have  dimension
$d_{j_0}=n+r$ for 
$s^{\prime\prime}$ near 0 in $\bR^q$. Using the
parametrizations
\thetag{2.2.10},
\thetag{2.2.13} and writing
$\Lambda=(z,\Lambda^\prime)$, we can express 
$N_{j_0}(s^{\prime\prime})$ in the coordinates
$(z,w)$ by
$$
w=v^{j_0}(z,\Lambda^\prime;s^{\prime\prime}),\tag3.4.3
$$ where
$\Lambda^\prime\in\bC^{(j-1)n}$. Since the defining
equations
\thetag{3.4.2} of $M$ at
$(0,0,s^{\prime\prime})$ depend algebraically on
$s^{\prime\prime}$, it follows that
$v^{j_0} (\cdot;s^{\prime\prime})$ also does (cf.
\thetag{2.2.10}--\thetag{2.2.12}  and
\thetag{1.2.13}--\thetag{2.2.15}). At a point
$(z,\Lambda^\prime;0)$ where 
$
\frac{\partial
v^{j_0}}{\partial\Lambda^\prime}
$ has maximal rank
$r=d_{j_0}-n$, we may assume  (by a change of
coordinates in the
$\Lambda^\prime$ space if necessary) that
$\Lambda^\prime=(\Lambda_1,
\Lambda_2)\in\bC^r\times\bC^{(j-1)n-r}$,
$v^{j_0}$ is independent of 
$\Lambda_2$, and
$
\frac{\partial v^{j_0}}{\partial\Lambda_1}
$ has rank $r$. Since the tangent plane of $X_0$ at
0 equals
$\{w^{\prime\prime}=0\}$, it follows from the
implicit function theorem that we can solve for
$\Lambda_1$ in the first $r$ equations of
\thetag{3.4.3}. We then substitute this into  the
last $q$ equations and find that we can express 
$N_{j_0}(s^{\prime\prime})$, for
$s^{\prime\prime}$ close to 0, as a graph
$$
w^{\prime\prime}=f(z,w^\prime;s^{\prime\prime})
\tag3.4.4
$$ near some point
$(z^1,{w^\prime}^1,f(z^1,{w^\prime}^1;
s^{\prime\prime}))$
with
$f(z,w^\prime;s^{\prime\prime})$ holomorphic
algebraic  in a neighborhood of $(z^1,w^1;0)$. Now,
since all the CR orbits near $p_0$ have the same
dimension, it follows from the Frobenius theorem that
they form a real analytic foliation of a neighborhood of
$p_0$ in $M$ (as we have noted before, Frobenius does
not imply that the orbits form a real {\it algebraic}
foliation even though the vector fields are algebraic).
Thus, there are $q$ real-valued, real analytic functions
$k=(k_1,...,k_q)$ on $M$ with linearly independent
differentials near $p_0$ such that every local CR
orbit near this point is of the form $\{(z,w)\in
M\:k(z,w)=c\}$ for some small
$c\in\bR^q$ (we may assume that
$k(0)=0$). Since
$(0,0,s^{\prime\prime})\in W_ {s^{\prime\prime}}$, we
have 
$$
W_{s^{\prime\prime}}=\{(z,w^\prime,w^{\prime\prime}))\in
M\:k(z,w^\prime,
w^{\prime\prime})=k(0,0,s^{\prime\prime})\}.
\tag3.4.5
$$ Clearly, these functions are CR and,
hence, they extend, near $0$ in $\bC^N$,  as
holomorphic functions which we again denote by
$k$. It follows that each
$X_{s^{\prime\prime}}$, for real $s^{\prime\prime}$
close to 0, is given by
$$
X_{s^{\prime\prime}}=\{(z,w^\prime,w^{\prime\prime}))\in
\bC^N\:k(z,w^\prime,
w^{\prime\prime})=k(0,0,s^{\prime\prime})\}.
\tag3.4.6
$$ Since the tangent plane of $X_0$ at 0 equals
$\{w^{\prime\prime}=0\}$, it follows that there is a
holomorphic function
$g(z,w^{\prime},s^{\prime\prime})$ near
$(0,0,0)$  with
$ g(0,0,s^{\prime\prime})\equiv
s^{\prime\prime}$ such that
$X_{s^{\prime\prime}}$,  for real
$s^{\prime\prime}$ close to 0, is given by
$$
w^{\prime\prime}=g(z,w^\prime,s^{\prime\prime}).
\tag3.4.7
$$ The maximal Segre set
$N_{j_0}(s^{\prime\prime})$ coincides with $X_
{s^{\prime\prime}}$ on a dense open subset of the
latter. Consequently, the algebraic representation
\thetag{3.4.4} of
$N_{j_0}(s^{\prime\prime})$, which is valid near the
point
$(z^1,{w^\prime}^1,f(z^1,{w^\prime}^1;
s^{\prime\prime}))$, implies that the holomorphic
function
$g(z,w^\prime, s^{\prime\prime})$ in (3.4.7) is in
fact algebraic. (The point
$(z^1,{w^{\prime}}^1)$  can be taken arbitrarily
close to 0.) Hence the algebraic
function
$f(z,w^\prime;s^{\prime\prime})$ can be continued to
an algebraic holomorphic function near
$(0,0;0)$. 

Now, as we noted above, we have the identity
$ f(0,0;s^{\prime\prime})\equiv
s^{\prime\prime}$ and hence 
$$
\frac{\partial f}{\partial
s^{\prime\prime}}(0,0;0)=I.\tag3.4.8
$$   Hence, we
may solve the equation
$$
w^{\prime\prime}=f(z,w^{\prime};
s^{\prime\prime})\tag3.4.9
$$ for $s^{\prime\prime}$ near the base point
$(z,w^\prime,w^{\prime\prime},
s^{\prime\prime})=(0,0,0,0)$. We obtain a 
$\bC^q$-valued algebraic function 
$h(z,w^\prime,w^{\prime\prime})$,
holomorphic near 
$(0,0,0)$, satisfying
$$
w^{\prime\prime}\equiv
f(z,w^{\prime};h(z,w^\prime,w^{\prime\prime}))
\tag3.4.10
$$ with
$h(0,0,s^{\prime\prime})=s^{\prime\prime}$. It
follows that the restriction of
$h(z,w^\prime,w^{\prime\prime})$ to 
$X_{s^{\prime\prime}}$ is constant and equals
$s^{\prime\prime}$. In  particular, since the CR
orbits
$W_{s^{\prime\prime}}=M\cap X_ {s^{\prime\prime}}$
(for
$s^{\prime\prime}\in\bR^q$ close to 0) cover a 
neighborhood of 0  in
$M$, the restriction of $h$ to $M$ is valued in
$\bR^q$. Indeed, we have
$h|_M=s^{\prime\prime}$ and, as a consequence, we
also have
$$
\partial h_1(0)\wedge...\wedge\partial h_q(0)\neq
0.\tag3.4.11
$$
The proof of Lemma 3.4.1 is complete. \qed\enddemo

\subhead 3.5. An example \endsubhead Consider the
 five 
dimensional real  algebraic submanifold $M
\subset\bC^4$ defined by
$$
\re Z_3 =0, \ \ \ \im Z_3 =|Z_1|^2, \ \ \ \im
Z_4=|Z_2|^2.
\tag3.5.1
$$ 
On the set
$\{(0,Z_2,0,X_4+i|Z_2|^2)\:Z_2\in\bC\,,\,X_4\in\bR\}$, 
$M$ is neither  generic nor CR, but outside
that set $M$ is  generic and 
holomorphically nondegenerate. The
function
$h_1(Z)=-iZ_3$ is real on $M$, but 
$M\cap\{Z\:h_1(Z)=c\}$, for real $c>0$, is
not minimal anywhere. Indeed,
$M\cap\{Z\:h_1(Z)=c\}$ is given by
$$
 |Z_1|^2 =c,\ \ Z_3= ic,\  \ \im Z_4 =|Z_2|^2,
\tag3.5.2
$$ which is not minimal since it is a product of a
circle and a three dimensional  surface. We leave
it to the reader to check that there is no germ at
$0$ of an algebraic holomorphic function $h$ which
is real on
$M$ and such that $\partial h(0)\wedge
\partial h_1(0)\not= 0$. 

Hence, we cannot apply Theorem 3.1.2 with
$p_0=0$. However, a straightforward
calculation reveals that the function 
$$ h_2(Z)=\frac{Z_1^2-iZ_3}{2Z_1}.\tag3.5.3
$$ is real on $M$, since
$h_2(Z)|_M=\re Z_1$. Near any point
$p_1=(ir,0,ir^2,0)\in M$, with
$r\in \bR$, the leaves 
$\{Z\:h_1(Z)=c_1\,,\,h_2(Z)=c_2\}$,
for $c=(c_1,c_2)\in\bC^2$ close  to $(r^2,0)$, are
equal to
$$
\{Z\:Z_1=c_2+\sqrt{c_2^2-c_1}\, , \, Z_3=ic_1\},
\tag3.5.4
$$ where the square root is chosen so that
$\sqrt{-1}=i$.  Assume now that  there is a
holomorphic map
$H\:\bC^4\mapsto\bC^{N^\prime}$ near 0,  generically
of rank 4, such that $H(M)$ is contained in a 5
dimensional real algebraic subset of
$\bC^{N^\prime}$. If we choose the point $p_1$ as
above with $r\not= 0$ to be in the domain of
definition of
$H$  then we may apply Theorem 3.1.2 in a
neighborhood of
$p_1$ since
$M\cap\{Z\:h_1(Z)=c_1\,,\,h_2(Z)=c_2\}$ is
minimal for
$Z$  near $p_1$,  and
$c \in \bR^2$  near $(r^2,0)$. Theorem 3.1.2
implies that $H$ is algebraic on the leaves 
$$
\{Z\:h_1(Z)=c_1\,,\,h_2(Z)=c_2\},\tag3.5.5
$$ which are the same as the leaves defined
by (3.5.4). More precisely, the proof of Theorem
3.1.2 implies that there are polynomials
$P_l(Z_2,Z_4,X; Z_1,Z_3)$ in
$(Z_2,Z_4,X)\in\bC^3$ with coefficients that are
holomorphic functions of $(Z_1,Z_3)$ near
$(ir,ir^2)$ such that (with
$H=(H_1,...,H_{N^\prime})$)
$$
P_l(Z_2,Z_4,H_l(Z);Z_1,Z_3)\equiv0\tag3.5.6
$$ holds for $Z$ near
$(ir,0,ir^2,0)$, for
$l=1,...,N^\prime$.  Since $H$ is holomorphic in a
neighborhood of 0, we can now apply Lemma 3.2.1
to conclude that $H$ is algebraic on  the leaves
$\{Z\:Z_1=Z_1^0\,,\,Z_3=Z_3^0\}$ for all
$(Z_1^0,Z_3^0)$ in a neighborhood of $(0,0)$. 
Note that as mentioned above, we could not apply
Theorem 3.1.2, as it is formulated, directly to this
example at
$p_0=0$. 

It should be noted that there exists 
a nonalgebraic mapping $H$ which is holomorphic
outside
$\{Z_1=0\}$, maps $M$ into itself, has generically
full rank, and which is algebraic on the
leaves $\{Z\:Z_1=Z_1^0\,,\,Z_3=Z_3^0\}$,  Indeed,
we may take
$$
H(Z_1,Z_2,Z_3,Z_4)=(e^{ih_2(Z_1,Z_3)-iZ_3/2}Z_1,Z_2,e^{-iZ_3}Z_3,Z_4).
\tag3.5.14
$$

\subhead 3.6.  Proofs of Theorems 1 through 4
\endsubhead

We begin by proving Theorem 1.  Condition (2) of 
Theorem 1 implies
that there is a point $p$ in $\acr$ at which $\acr$ is
generic; by Lemma 3.4.1 we may also assume that
$\acr$ is minimal at $p$.  The proof of Theorem 1 then
follows from Corollary 3.1.4.

For Theorem 4, we note first that Theorem 2.2.1 states
that the CR orbits and their intrinsic complexifications
are all algebraic. The rest of the proof of the theorem
follows from Theorem 1, since any biholomorphism must
map a CR orbit onto a CR orbit.

Now we shall prove Theorems 2 and 3. By Proposition 1.4.1,
condition (1) of Theorem 1 is equivalent to
condition (i) of Theorem 3.  We  first show
that holomorphic degeneracy implies property (3)
of Theorem 2.   

\proclaim {Proposition 3.6.1} Let
$A$ be an irreducible real algebraic subset of
$\cnn$. If either {\rm (i)} or {\rm (ii)} of Theorem 3
does not hold, then {\rm (3)} of Theorem 2 holds.
\endproclaim
\demo {Proof} Assume first that (i) does not
hold and let $\po \in
\acr$.    By Proposition 1.4.1, the definition
of holomorphic degeneracy, and the
observations in the proof of Proposition
1.4.1, there exists a nontrivial holomorphic
vector field
$X$ of the form (1.4.1) tangent to
$A$ with coefficients algebraic holomorphic
near $\po$.  Without loss of generality, we
may assume
$X(\po) = 0$.  The proof now is essentially
the same as that of the hypersurface case
([BR3, Proposition 3.5]). We take the
complex flow of the vector field
$X$ or, if necessary, of $fX$, where $f$ is a
germ of a nonalgebraic holomorphic function
at $\po$ to find the desired germ of
biholomorphism satisfying (3). See [BR3] for
details. 

Assume now that (ii) does not hold, and let
$\po\in \acr$.  Since $A$ is not generic at
$\po$, there exists an algebraic holomorphic
proper submanifold in
$\cnn$ containing $\acr$.  After an algebraic
holomorphic change of coordinates, we may
assume that
$\po = 0$ and that $A$ is contained in the
complex hyperplane $Z_N = 0$ near
$0$.  To prove that (3) holds, it suffices to
take the mapping
$H_j(Z) = Z_j,\ j =1, \ldots, N-1$, and
$H_N(Z)= Z_Ne^{Z_N}$. This proves Proposition
3.6.1.  
  $\square$ 
\enddemo

We now prove the last statement of Theorem 3.
A homogeneous submanifold
$M$ of
$\cnn$ of codimension $d$ is given by
$$ M = \{Z \in \cnn: \r_j(Z,\overline Z) = 0, j
=1,\ldots,d\}, \tag 3.6.1
$$  where the $\r_j$ are real valued
polynomials weighted homogeneous with respect to the
weights $\nu_1\le\ldots\le\nu_N$ (see \S 2.3). Let 
 $r_1\le
\ldots\le r_d$ be the degrees of homogeneity of
the polynomials $\rho_1, \ldots, \rho_d$, i.e., for
$t > 0$
$$\r_j(t^{\nu_1}Z_1,\ldots,t^{\nu_N}Z_N) =
t^{r_j}\r_j(Z,\overline Z), \ \ j=1,\ldots,d.
\tag 3.6.2$$ We also assume that 
$$ d\r_1(0) \wedge \ldots
\wedge d\r_d(0) \not= 0. \tag 3.6.3
$$

\proclaim {Lemma 3.6.2} Let $M$ be a
homogeneous generic submanifold of
$\cnn$  which is not minimal at
$0$.  Then there exists a holomorphic
polynomial $h$ in
$\cnn$, with $h|_M$ nonconstant and real
valued.  
\endproclaim

\demo {Proof} The homogeneous manifold $M$ is
generic (at $0$ and hence at all points) if,
in addition,  to (3.6.3) we have
$$ \pa\r_1(0) \wedge
\ldots
\wedge \pa\r_d(0) \not= 0. \tag 3.6.4
$$ The reader can easily check that if
$M$ is a generic homogeneous manifold of
codimension $d$, after a linear holmorphic
change of coordinates 
$Z = (z,w)$,  $M$ can be written in the form
$$w=Q(z,\overline z,\overline w), \
\hbox {\rm with} \ Q_j(z,\overline z,\overline w) = \overline w_j +
q_j(z,\overline z,\overline w_1,\ldots,
\overline w_{j-1}), \tag 3.6.5$$
$j=1,\ldots, d$, with
$q_j$  a weighted homogeneous polynomial of
weight
$r_j$. Here $Q$ is complex valued and
satisfies (1.1.5). After a further weighted
homogeneous change of holomorphic
coordinates, we may assume that the
coordinates $(z,w)$ are normal, i.e. (1.1.4)
holds.

 As in \S 2, we let
$M^k$ be the projection of $M$ in
$\bC^{n+k-1}$, $k = 2,\ldots,d+1$.  Each $M^k$ is
defined by the first $k-1$ equations in
(3.6.5). If  the hypersurface
$M^2
\subset
\bC^{n+1}$ is not minimal at $0$, then
necessarily $q_1(z,\overline z) \equiv 0$, and we may
take $h(z,w) = w_1$.  If not, we let $\ell \le
d $ be the smallest integer for which
$M^{\ell}$ is minimal at $0$, but
$M^{\ell+1}$ is not minimal at $0$. Then the
CR orbit $W$ of $0$ in the generic manifold
$M^{\ell+1}$ is a proper CR submanifold of
$M^{\ell+1}$ of CR dimension $n$. It must be a
holomorphic graph over $M^\ell$ in
$\bC^{n+\ell}$.  That is, $W$ is given by
(3.6.5) for $1 \le j \le \ell-1$ and
$w_{\ell+1} = f(z,w_1,\ldots ,w_{\ell-1})$.
Since
$W \subset M^{\ell+1}$ we must also have $\im
f(z,w_1,\ldots,w_{\ell-1})|_M =
(1/2i)q_{\ell}(z,\overline z,\overline w_1,\ldots
\overline w_{\ell-1})|_M $ .  The reader can check that
this implies that $f(z,w_1,\ldots,w_{\ell-1})$ is
independent of $z$ and is a weighted
homogeneous holomorphic polynomial, and the
function
$h(z,w) = w_{\ell} -
f(z,w_1,\ldots,w_{\ell-1})$ satisfies the
conclusion of the lemma.
  $\square$ 
\enddemo

\proclaim {Proposition 3.6.3}  Let $M$ be a
homogeneous generic submanifold of
$\cnn$ which is not minimal at
$0$. Then for any $\po \in M$, there exists a
nonalgebraic holomorphic map $H$ from $\cnn$
into itself with $H(\po) = \po$,
$H(M)\subset M$, with $\J H(\po)
\not= 0$.  
\endproclaim

\demo {Proof}  By Lemma 3.6.2, there exists a
nonconstant holomorphic polynomial 
$h$ with $h|_M$ real.  We may also assume
$h(\po) = 0.$ The reader can easily check
that  the map defined by
$$H_j(Z) = e^{\nu_j h(Z)}Z_j, \ \ j=1,\ldots,
N $$ satisfies the desired conclusion of the
proposition.
$\square$ 
\enddemo

\enddocument
\end